\documentclass[onefignum,onetabnum]{siamonline220329}
\usepackage{lipsum}
\usepackage{amsfonts}
\usepackage{graphicx}
\usepackage{epstopdf}
\usepackage{algorithmic}
\ifpdf
  \DeclareGraphicsExtensions{.eps,.pdf,.png,.jpg}
\else
  \DeclareGraphicsExtensions{.eps}
\fi
\usepackage{tikz}
\usepackage{pgfplots}
\pgfplotsset{compat=1.15}
\usetikzlibrary{intersections, pgfplots.fillbetween}
\usepackage[numbers,sort&compress]{natbib}

\usepackage{enumitem}
\setlist[enumerate]{leftmargin=.5in}
\setlist[itemize]{leftmargin=.5in}

\newsiamremark{remark}{Remark}
\newsiamremark{example}{Example}
	\crefname{example}{Example}{Examples}
\newsiamremark{hypothesis}{Hypothesis}
	\crefname{hypothesis}{Hypothesis}{Hypotheses}
\newsiamremark{assumption}{Assumption}
	\crefname{assumption}{Assumption}{Assumptions}
\newsiamremark{property}{Property}
	\crefname{property}{Property}{Properties}
\newsiamthm{claim}{Claim}

\headers{Gradient Descent with Diminishing Step Sizes
}{Patel \& Berahas}

\title{Gradient Descent in the Absence of Global Lipschitz Continuity of the Gradients\thanks{Preprint
\funding{This work is supported by the Wisconsin Alumni Research Foundation and the Office of Naval Research.}}
}

\author{Vivak Patel\thanks{Department of Statistics, University of Wisconsin -- Madison, Madison, WI, USA 
  (\email{vivak.patel@wisc.edu}, \url{vivakpatel.org}).}
\and Albert S.  Berahas\thanks{Department of Industrial and Operations Engineering, University of Michigan -- Ann Arbor, Ann Arbor, MI USA 
  (\email{albertberahas@gmail.com}, \url{aberahas.engin.umich.edu).}}
}

\usepackage{amsopn}
\usepackage{amssymb}
\usepackage{tcolorbox}
\tcbuselibrary{breakable}
\tcbset{%any default parameters
  width=0.99\textwidth,
  halign=justify,
  center,
  breakable,
  colback=white    
}
\usepackage{booktabs}

\newcommand{\norm}[1]{\left\Vert #1 \right\Vert}
\newcommand{\inlinenorm}[1]{\Vert #1 \Vert}

\newcommand{\flb}{F_{l.b.}}

\newcommand{\cmp}[2]{#1^{(#2)}}

\ifpdf
\hypersetup{
  pdftitle={Gradient Descent with Diminishing Step Sizes},
  pdfauthor={Patel \& Berahas}
}
\fi

\begin{document}

\maketitle

% REQUIRED
\begin{abstract}
Gradient descent (GD) is a collection of continuous optimization methods that have achieved immeasurable success in practice. Owing to data science applications, GD with diminishing step sizes has become a prominent variant. While this variant of GD has been well-studied in the literature for objectives with globally Lipschitz continuous gradients or by requiring bounded iterates, objectives from data science problems do not satisfy such assumptions. Thus, in this work, we provide a novel global convergence analysis of GD with diminishing step sizes for differentiable nonconvex functions whose gradients are only locally Lipschitz continuous. Through our analysis, we generalize what is known about gradient descent with diminishing step sizes including interesting topological facts; and we elucidate the varied behaviors that can occur in the previously overlooked divergence regime. Thus, we provide the most general global convergence analysis of GD with diminishing step sizes under realistic conditions for data science problems.
\end{abstract}

% REQUIRED

% REQUIRED

\section{Introduction} \label{section-introduction}
Proposed nearly two centuries ago \cite{cauchy1847,hadamard1908,curry1944,lemarechal2012}, gradient descent is a set of canonical continuous optimization methods that have achieved immeasurable success in a plethora of applications (e.g.,~\cite{dommel1968,bottou2010,iyengar2013}). Owing to their prominence and utility in data science, gradient descent methods have continued to grow in variety, and their theory has received renewed interest by the optimization and data science communities for problems in this area (e.g., \cite{du2018,du2019,lee2019,oymak2019}). In particular, gradient descent with pre-scheduled step sizes has become popular owing to the additional expense of using line search techniques for data science problems. Correspondingly, the theory of gradient descent with pre-scheduled step sizes has grown in a number of interesting directions including new local convergence rate analyses (e.g., \cite{karimi2016,du2019}) and saddle-point avoidance analyses (e.g., \cite{lee2016,du2017,jin2018}).

That said, the more fundamental global convergence analysis of gradient descent with pre-scheduled step sizes has lagged owing to two challenges. First, gradient descent with pre-scheduled step sizes does not guarantee a monotonic reduction in the objective function (c.f., Armijo's Method \cite{armijo1966}), which is the key ingredient used to analyze such methods via Zoutendjik's approach \cite{zoutendijk1960}. 
Second, because of the nonconvexity of common data science problems,\footnote{Canonical data science problems such as Poisson regression, linear three-or-more layer feed forward networks, and linear three-or-more time horizon recurrent networks fail to possess globally Lipschitz continuous gradients or uniformly continuous gradients when trained using standard loss functions \citep[\S 1]{patel2021global}.}  
the analysis of gradient descent cannot leverage uniform continuity of the gradient function, global Lipschitz continuity of the gradient function or presuppose that its iterates remain in a bounded region for a function with a locally Lipschitz continuous gradient,\footnote{If the iterates remain in a bounded region, then compactness and the local Lipschitz continuity condition would imply a global Lipschitz continuous constant in the bounded region.} which are instrumental assumptions for overcoming the previous challenge \cite{reddi2016,reddi2016stochastic}. As a result of the latter challenge, typical analysis approaches for global convergence of gradient descent fall short (see \cref{subsection-literature-results} for an overview).
In fact, even the new vogue for analysis in machine learning, \textit{the continuous approach} \cite{ljung1977,benaim1999,benveniste2012,mertikopoulos2020}, falls short, as this approach requires compactness of the image space of the iterates in \cite[Theorem 3.2]{benaim1999}, boundedness of iterates \cite[Theorem 2]{fort1996}, or global Lipschitz continuity of the gradient of the objective function \cite[Assumption 1]{mertikopoulos2020}.
To summarize, to the best of our knowledge, existing global convergence analyses of gradient descent with diminishing step sizes do not apply to canonical, nonconvex, differentiable data science problems. 

To address this shortcoming, we generalize recently developed techniques for the analysis of stochastic gradient descent \cite{patel2020,patel2021,patel2021global} to analyze gradient descent with diminishing step sizes for nonconvex optimization problems that are bounded from below and whose gradient is \textit{locally Lipschitz continuous}, which are more realistic assumptions for canonical data science problems \cite[\S A]{patel2021global}. Our analysis has several important contributions.
\begin{remunerate}
\item First, we present a novel upper-bound model, which can be used under milder assumptions that are appropriate for data science problems (see \cref{subsection-literature-results} for a discussion, and \cref{result-recursion-gap} for the result). This upper-bound model is directly useful in analyzing many other algorithms for unconstrained optimization, and the strategies used to prove the result seem useful for analyzing algorithms for constrained optimization.
\item Second, our analysis provides counterexamples to what is known about gradient descent with diminishing step sizes. Specifically, previous results (e.g., \citep[Proposition 1.2.4]{bertsekas2016}) showed that, under a global Lipschitz continuity assumption on the gradient, the iterates tend to a region where the gradient is zero; the objective function converges to a finite limit; and, if the iterates remain bounded, then the iterates converge to a stationary point. Our analysis, under the more realistic local Lipschitz continuity assumption on the gradient, offers a correction to this view---that the gradient function can remain bounded away from zero and the objective function can diverge (see explicitly constructed examples in \cref{section-divergence-example}). 
\item Our analysis addresses a preliminary question about gradient descent and nonconvexity: \emph{given a relatively arbitrary objective function, can its nonconvexity cause gradient descent to behave erratically in a region?} Despite the general nonconvexity allowed by our assumptions, we show that the limit supremum and limit infimum of the objective function evaluated at the iterates must tend to each other if the iterates remain in a region for long enough---even if they eventually escape; we also show that the limit of the gradient function evaluated at the iterates must tend to zero if the iterates remain in a region for long enough---even if they eventually escape---(see \cref{result-asymptotics-of-iterates}). A more interesting question is whether such a statement holds uniformly over important subsets of nonconvex objective functions of the ones considered here.\footnote{Such functions must be beyond those which have globally Lipschitz continuous gradients, and still be valid for data science problems. To the best of our knowledge, the appropriate notion of continuity of the gradients has not been determined.} Our analysis at least gives hope that such a statement may be true.
\item Our analysis adds several topological insights to what is known (e.g., \citep[Proposition 1.2.4]{bertsekas2016}). Primarily, we show, the subsequential limit points of the iterates are a connected set that is either a singleton or infinite. Moreover, if the set is infinite, we conclude that it cannot contain an open set.
\end{remunerate}

\begin{tcolorbox}
Thus, to the best of our knowledge, our results provide the most general and complete global convergence analysis of gradient descent with diminishing step sizes under realistic assumptions for nonconvex, differentiable optimization problems that arise in data science.
\end{tcolorbox}

The remainder of this work is organized as follows. In \cref{section-gradient-descent}, we specify the class of nonconvex optimization problems of interest and the precise form of gradient descent with diminishing step sizes. In \cref{section-no-summable}, we analyze the behavior of gradient descent with diminishing step sizes. In \cref{section-divergence-example}, we construct examples that elucidate the possible behaviors of gradient descent with diminishing step sizes in the divergence regime. 
Final remarks are given in \cref{section-conclusion}.

\section{Gradient Descent} \label{section-gradient-descent}
We begin by introducing the general class of optimization problems that we consider in this work. Then, we specify the precise form of gradient descent with diminishing step sizes. With the problem class and procedure specified, we describe relevant analysis approaches in the literature.

\subsection{Optimization Problem}

To cover a variety of canonical problems in data science \cite[\S A]{patel2021global}, consider the optimization problem 
\begin{equation} \label{eqn-problem}
\min_{x \in \mathbb{R}^p} F(x),
\end{equation}
under the following assumptions.
\begin{assumption} \label{assumption-lower-bound}
The objective function, $F:\mathbb{R}^p \to \mathbb{R}$, is bounded from below by a constant $\flb$.
\end{assumption}
\begin{assumption} \label{assumption-gradient}
The gradient function, $\dot F(x) = \nabla F(z) \vert_{z = x}$, exists $\forall x \in \mathbb{R}^p$, and is locally Lipschitz continuous. 
\end{assumption}

For our context, we use the following definition of local Lipschitz continuity.
\begin{definition} \label{definition-local-lipschitz}
A function $G: \mathbb{R}^p \to \mathbb{R}^p$ is locally Lipschitz continuous if for every $x \in \mathbb{R}^p$ there exists an open ball of $x$, $\mathcal{N}$, and there exists $L \geq 0$, such that, $\forall y,z \in \mathcal{N}$,
\begin{equation} \label{eqn-lipschitz-continuity}
\frac{ \norm{ G(y) - G(z) }_2 }{ \norm{ y - z}_2 } \leq L.
\end{equation}
\end{definition}
Equivalently, $G$ is locally Lipschitz continuous if for every compact set $\mathcal{C} \subset \mathbb{R}^p$, there exists $L \geq 0$ such that \cref{eqn-lipschitz-continuity} holds for all $y,z \in \mathcal{C}$. This well-known statement is shown in \cref{result-equivalent-lipschitz-definitions}.

To give an example of the broad applicability of \cref{assumption-gradient}, any optimization problem whose objective function is twice continuously differentiable immediately satisfies \cref{assumption-gradient}. This well-known statement is given formally in \cref{result-hessian-implies-lipschitz}.

\subsection{Gradient Descent with Diminishing Step Sizes}

Now, suppose we apply gradient descent with diminishing step-sizes to solve \cref{eqn-problem}. Specifically, given $x_0 \in \mathbb{R}^p$, we generate a sequence $\lbrace x_k : k \in \mathbb{N} \rbrace$ according to
\begin{equation} \label{eqn-recursion-iterate}
x_{k+1} = x_k - M_k \dot F(x_k),
\end{equation}
where $M_k$ satisfies some of the following properties:
\begin{property} \label{property-spd}
$\lbrace M_k : k + 1 \in \mathbb{N} \rbrace \subset \mathbb{R}^{p \times p}$ are symmetric positive definite matrices;
\end{property}
\begin{property} \label{property-diverge}
$\sum_{k=0}^\infty \lambda_{\min}(M_k)$ diverges, where $\lambda_{\min}(M_k)$ denotes the smallest eigenvalue of $M_k$; and
\end{property}
\begin{property} \label{property-diminishing}
$\lim_{k \to \infty} \lambda_{\max}(M_k) = 0$ where $\lambda_{\max}(M_k)$ denotes the largest eigenvalue of $M_k$.
\end{property} 

\cref{property-spd,property-diverge,property-diminishing} are a matrix-valued generalization of classical diminishing step size requirements \cite[Proposition 1.2.4]{bertsekas2016}. Moreover, \cref{property-spd,property-diverge,property-diminishing} are enough to show that the objective function evaluated at the iterates converges, and to show that the limit infimum of the norm of the gradient function evaluated at the iterates converges to zero (see \Cref{result-zoutendjik}). To show that the gradient function converges to zero, these properties will be augmented with the following:
\begin{property} \label{property-conditioning}
There exists $\kappa \geq 1$ such that  $\lambda_{\max}(M_k) / \lambda_{\min}(M_k) \leq \kappa$ for all $k + 1 \in \mathbb{N}$.
\end{property}

Of interest, \Cref{property-spd,property-diverge,property-diminishing} can potentially account for adaptive step size selection procedures that exist in the literature, namely those that do not make use of objective function information. For example, \cref{property-spd,property-diverge,property-diminishing} can apply to the method of \cite{bertsekas2016} (with $\lambda=1$), which combines incremental nonlinear least squares, the Gauss-Newton method, and the Extended Kalman Filter. However, especially in the nonlinear case, \cref{property-diverge} would be difficult to verify without assuming something akin to what is called persistent excitation in the control literature \citep{johnstone1982exponential,bittanti1990convergence,parkum1992recursive,cao2003exponential}. 
Indeed, in the objective-free first-order optimization (e.g., AdaGrad-type methods), this persistent excitation condition often manifests through a combination of assumptions about the optimization problem (e.g., bounded gradients) and the diagonal or identity-scaling choice of $\lbrace M_k : k+1 \in \mathbb{N} \rbrace$ \citep{wu2018wngrad,ward2020adagrad,defossez2020simple,grapiglia2022adaptive,gratton2022convergence,gratton2022first,gratton2022parametric}.

\subsection{Important Analysis Approaches in the Literature} \label{subsection-literature-results}

With the problem and algorithm established, we briefly review two important analysis frameworks in the literature with respect to simple objective functions satisfying \cref{assumption-lower-bound,assumption-gradient}: $|x|^3$ and $\exp(x)$. Note, these two examples are essential components in verifying that canonical data science problems neither have globally Lipschitz continuous gradients nor uniformly continuous gradients \citep[see][\S A]{patel2021global}.

In one analysis framework for trust region methods \citep[e.g.,][]{more1983recent,ke1998class}, continuity of the gradient function, properties of the algorithm, and evaluations of the objective function are needed to show that the limit infimum of the gradient function evaluated at the iterates goes to zero. Furthermore, assuming uniform continuity of the gradient function allows for the conclusion: the limit of the gradient function evaluated at the iterates goes to zero \citep{more1983recent,ke1998class,zhang2006superlinearly,zhang2007trust,cartis2011adaptiveI,cartis2011adaptiveII}.\footnote{In \citep[AF.2]{cartis2011adaptiveI}, uniform continuity implies the needed property.} While continuity of the gradient function certainly holds for our two example objectives, neither of them satisfy uniform continuity of the gradient function. Moreover, in our context, gradient function information is not combined with objective function information to ensure sufficient decay at each step, which limits our ability to use the assumption of continuity of the gradient in place of \cref{assumption-gradient}.\footnote{This raises the question of how much objective function information is really needed in order to ensure similar results as trust-region without substantially increasing computational costs for data science problems.}

In the other analysis framework espoused by \citep[Proposition 1.2.4]{bertsekas2016}, \citep[Theorem 3.2]{nocedal2006}, and \citep[Lemma 10.4]{beck2017first}, the essential ingredient is a global upper-bound model for the objective function,
\begin{equation} \label{equation-global-upper-bound-model}
F(y) \leq F(x) + \dot F(x)^\intercal(y - x) + \frac{L}{2} \norm{ y - x}_2^2, ~\forall y, x \in \mathbb{R}^p,
\end{equation}
where $L$ is a fixed constant that arises from the assumption that the gradient function is \textit{globally Lipschitz continuous} (i.e., $L$ is the same regardless of $x \in \mathbb{R}^p $ and $\mathcal{N}$ in \cref{definition-local-lipschitz}). Indeed, this global upper-bound model is commonly used in recent analyses, both deterministic and stochastic \citep{wu2018wngrad,curtis2019stochastic,curtis2020adaptive,ward2020adagrad,defossez2020simple,grapiglia2022adaptive,gratton2022convergence,gratton2022first,gratton2022parametric}. 
This global upper-bound model is actualized by replacing $y$ with $x_{k+1}$, $x$ with $x_k$, and rewriting the right-hand-side strictly in terms of quantities depending on $x_k$. Then, the upper-bound model is manipulated to show that the objective function is decreasing. Unfortunately, such a global upper-bound model does not apply to the two simple example objective functions, which renders such analyses inapplicable to common data science problems.

In \citep[Exercise 1.2.5]{bertsekas2016}, this global upper-bound model is relaxed to the case where such an $L$ exists for every level set of the objective function and assumes every level set is bounded. In this case, this relaxed upper-bound model can then be used to establish that, if a gradient descent procedure remains in a level set, then the objective function converges to a finite value and the gradient function converges to zero. Indeed, this relaxed upper-bound model can account for $|x|^3$, but it cannot account for $\exp(x)$ nor our example in \cref{section-divergence-example}, which has bounded level sets yet the iterates never remain in any level set. Hence, even this relaxation cannot account for the types of problems that satisfy \cref{assumption-lower-bound,assumption-gradient}. 

Our approach can be viewed as a generalization of \citep[Exercise 1.2.5]{bertsekas2016} as we can use \cref{assumption-gradient} to write a valid upper-bound model for \textit{any two points in $\mathbb{R}^p$, even though we only assume local Lipschitz continuity of the gradient} (see \cref{result-recursion-gap} and \cref{example-limitations-of-classical-upper-bound-model}).
We now introduce this analysis approach.

\section{Global Convergence Analysis} \label{section-no-summable}
Here, we study the global convergence of gradient descent, \cref{eqn-recursion-iterate}, with diminishing step sizes satisfying \cref{property-spd,property-diverge,property-diminishing} on a general class of nonconvex functions as defined by \cref{assumption-lower-bound,assumption-gradient}. Our main conclusion is that, despite the allowed nonconvexity of a problem, the objective function and gradient function at the iterates are either stabilizing or the iterates must continually tend further away from the origin. Thus, if we somehow know that the iterates remain bounded, then they must converge to a stationary point.

To prove these claims, our main innovation is to analyze the gradient descent procedure under a stopping time framework, which is a theoretical construction that allows us to analyze the without modifying it. We enumerate the steps in our analysis here.
\begin{remunerate}
\item In \cref{subsection-recursion-gap}, we establish a novel upper bound model based on stopping times to relate the optimality gaps of two arbitrary points even under local Lipschitz continuity of the gradient function (see \cref{result-recursion-gap}). We then simplify this statement when we substitute the two arbitrary points with consecutive iterates generated by the gradient descent procedure with diminishing step sizes (see \cref{result-recursion-gap-simplified}).
\item In \cref{subsection-zoutendjik}, we apply Zoutendjik's analysis approach \cite{zoutendijk1960}. We show that the limit supremum and limit infimum of the objective function evaluated at the iterates must tend to each other if the iterates remain in a region for long enough (even if they eventually escape). We also show that the limit infimum of the gradient function evaluated at the iterates must tend to zero if the iterates remain in a region for long enough (even if they eventually escape). 
\item In \cref{subsection-convergence-gradient}, we strengthen the preceding statement: we show that the limit of the gradient function evaluated at the iterates tends to zero if the iterates remain in a region for long enough (even if they eventually escape).
\item In \cref{subsection-topological-iterates}, we establish topological properties of the iterates when their subsequential limits are a bounded set. In particular, we establish the well-known results that the limit points of the iterates converge to a closed set where the gradient function is zero, and we establish---to the best of our knowledge---the novel result that this set is connected and cannot contain an open set (see \cref{result-asymptotics-of-iterates}). In other words, when it converges, gradient descent with diminishing step sizes either tends to a single point or an infinite set that must, in a sense, lack volume. Moreover, gradient descent with diminishing step sizes cannot have a cycle nor can it converge to a limit cycle with a finite number of points.
\end{remunerate}

We turn our attention to the divergence regime in \cref{section-divergence-example}.

\subsection{A Relationship for the Optimality Gap} \label{subsection-recursion-gap}

We now establish an upper-bound inequality for the optimality gap between two points in $\mathbb{R}^p$ under local Lipschitz continuity (see \cref{subsection-literature-results}). To establish this result, we make use of a technique from probability theory that analyzes stochastic processes under stopping times. For the deterministic equivalent, we define for an arbitrary point $x \in \mathbb{R}^p$ and $R \geq 0$,
\begin{equation}
\pi_x(R) = \begin{cases}
1 & \inlinenorm{x}_2 \leq R \\
0 & \text{otherwise}.
\end{cases}
\end{equation}

\begin{lemma} \label{result-recursion-gap}
Suppose $F:\mathbb{R}^p \to \mathbb{R}$ satisfies \cref{assumption-lower-bound,assumption-gradient}. Then, for all $ R \geq 0$ there exists a constant $C_R > 0$ such that, for all $x,y \in \mathbb{R}^p$,
\begin{equation}
[F(y) - \flb]\pi_y(R)\pi_x(R) \leq \left[F(x) - \flb - \dot F(x)^\intercal (y-x) + C_R \norm{ y-x}_2^2\right]\pi_x(R). 
\end{equation}
\end{lemma}

At first glance, we might think that \cref{result-recursion-gap} can be proved by combining \cref{equation-global-upper-bound-model} with $L \geq 0$ specific to the radius $R > 0$ of interest to show
\begin{equation} 
[F(y) - \flb] \pi_y(R)\pi_x(R) \leq \left[ F(x) - \flb - \dot F(x)^\intercal (y-x) + \frac{L}{2} \norm{y-x}_2^2 \right] \pi_y(R)\pi_x(R),
\end{equation}
and then using $\pi_y(R)\pi_x(R) \leq \pi_x(R)$ to upper bound the right hand side to conclude,
\begin{equation} \label{equation-wrong-upper-bound-model}
[F(y) - \flb] \pi_y(R)\pi_x(R) \leq \left[ F(x) - \flb - \dot F(x)^\intercal (y-x) + \frac{L}{2} \norm{y-x}_2^2 \right] \pi_x(R).
\end{equation}
Unfortunately, it is the last step that can be problematic as the the right hand side can become negative, which produces a false inequality. The following example illustrates the issue.

\begin{example} \label{example-limitations-of-classical-upper-bound-model}
Consider
\begin{equation}
F(x) = \begin{cases}
10(1-x) & x \leq 1 \\
\frac{1}{10x - 9} - 1 & x > 1,
\end{cases}
\quad\text{for which}\quad
\dot F(x) = \begin{cases}
-10 & x \leq 1 \\
\frac{-10}{(10x - 9)^2} & x > 1,
\end{cases}
\end{equation}
which is bounded from below and for which $\dot F(x)$ is globally Lipschitz continuous. If we now set $R = 1$, $x = 1$, then we see that $L = 0$ on $[-1,1]$ and $\pi_x(1) = 1$. If we now select $y=11$ (which would be the iterate generated by a gradient descent procedure at $x=1$ with step size $1$) then $\pi_y(1) = 0$. Plugging this into \cref{equation-wrong-upper-bound-model},
$0 = (F(11) + 1) 0 \leq ( 0 + 1 - 100 ) 1 = -99$,
which is false. Hence, proving \cref{result-recursion-gap} requires a little more care as we show below.
\end{example}

\begin{proof}
First, for any $R \geq 0$ define $L_R$ to be the Lipschitz constant for the gradient in the closed ball of radius $R$ around the point $0 \in \mathbb{R}^p$, which is well defined by \cref{assumption-gradient,result-equivalent-lipschitz-definitions}. For any fixed $\delta > 0$, it readily follows that $L_R \leq L_{R+\delta}$. Second, let $L(y,x)$ be the Lipschitz constant of the gradient in a closed ball of radius $\inlinenorm{ y - x}_2$ around the point $x$. Finally, for any $R \geq 0$, define $G_R$ to be the maximum  $\inlinenorm{ \dot F(x)}_2$ for all $x$ in a closed ball of radius $R$ around the point $0 \in \mathbb{R}^p$. Now, let $y,x \in \mathbb{R}^p$ be arbitrary.

By Taylor's remainder theorem,
\begin{equation}
\begin{aligned}
F(y) - \flb &= F(x) - \flb + \dot F(x)^\intercal (y - x) \\
&\quad + \int_0^1 \left[ \dot F(x + t(y-x)) - \dot F(x) \right]^\intercal (y - x) dt.
\end{aligned}
\end{equation}
By applying \cref{assumption-gradient} to the last term,
\begin{equation}\label{eqn-no-lipschitz-control}
    F(y) - \flb \leq F(x) - \flb + \dot F(x)^\intercal (y - x) + \frac{L(y,x)}{2} \norm{ y-x}_2^2\\
\end{equation}
Note, to understand why we must keep going at this point in the proof, see \cref{discussion-no-lipschitz-control}. We now introduce $\pi_y(R)$ and $\pi_x(R)$ into \cref{eqn-no-lipschitz-control}. That is,
\begin{equation}
\begin{aligned}
&[F(y) - \flb] \pi_y(R) \pi_x(R) \\
&\quad \leq \left[F(x) - \flb - \dot F(x)^\intercal (y-x) + \frac{L(y,x)}{2} \norm{y-x}_2^2\right] \pi_y(R) \pi_x(R).
\end{aligned}
\end{equation}
If $\pi_y(R)\pi_x(R) = 1$, then $\inlinenorm{y}_2 \leq R$ and $\inlinenorm{x}_2 \leq R$. Thus, $L(y,x) \leq L_R \leq L_{R+\delta}$. When $\pi_y(R)\pi_x(R) = 0$, then both sides are trivially zero. Therefore,
\begin{equation}
\begin{aligned}
&[F(y) - \flb] \pi_y(R)\pi_x(R) \\
& \quad \leq \left[F(x) - \flb - \dot F(x)^\intercal (y-x) + \frac{L_{R+\delta}}{2} \norm{y-x}_2^2\right] \pi_y(R)\pi_x(R).
\end{aligned}
\end{equation}

Now, we want $\pi_x(R)$ alone on the right hand side. So, we simply add and subtract a term involving $\pi_x(R)$, and study the difference term. That is,
\begin{equation} \label{eqn-add-subtract-event}
\begin{aligned}
&[F(y) - \flb] \pi_y(R)\pi_x(R) \\
&\quad \leq \left[F(x) - \flb - \dot F(x)^\intercal(y-x) + \frac{L_{R+\delta}}{2} \norm{y-x}_2^2\right] \pi_x(R) \\
&\qquad  + \left[F(x) - \flb - \dot F(x)^\intercal (y-x) + \frac{L_{R+\delta}}{2} \norm{y-x}_2^2\right] [\pi_y(R)\pi_x(R) - \pi_x(R)].
\end{aligned}
\end{equation}
We now have two cases to upper bound the last term of \cref{eqn-add-subtract-event}. Note, $\pi_y(R)\pi_x(R) - \pi_x(R) \leq 0$.
\paragraph{Case 1} If
\begin{equation}
F(x) - \flb - \dot F(x)^\intercal (y-x) + \frac{L_{R+\delta}}{2} \norm{y-x}_2^2 \geq 0,
\end{equation}
then
\begin{equation}
\left[F(x) - \flb - \dot F(x)^\intercal (y-x) + \frac{L_{R+\delta}}{2} \norm{y-x}_2^2\right][\pi_y(R)\pi_x(R) - \pi_x(R)] \leq 0.
\end{equation}
Hence, in this case, we can upper bound the last term in \cref{eqn-add-subtract-event} by any non-negative term.

\paragraph{Case 2} If
\begin{equation} \label{eqn-add-subtrct-event-case-2}
F(x) - \flb - \dot F(x)^\intercal (y-x) + \frac{L_{R+\delta}}{2} \norm{y-x}_2^2 < 0,
\end{equation}
then, using $\pi_y(R)\pi_x(R) \leq \pi_x(R)$,
\begin{equation}
\left[F(x) - \flb - \dot F(x)^\intercal(y-x) + \frac{L_{R+\delta}}{2} \norm{y-x}_2^2\right][\pi_y(R)\pi_x(R) - \pi_x(R)]\geq 0.
\end{equation}
Thus, we need only to find a lower bound for the first term in the product to upper bound the entire term. Specifically, 
\begin{equation}
- \norm{ \dot F(x)}_2 \norm{ y-x}_2 \leq F(x) - \flb - \dot F(x)^\intercal (y-x) + \frac{L_{R+\delta}}{2} \norm{y-x}_2^2.
\end{equation}
Now, when $\pi_y(R)\pi_x(R) < \pi_x(R)$, $\inlinenorm{x}_2 \leq R$. Moreover, if \cref{eqn-add-subtrct-event-case-2} holds, then $R+\delta < \inlinenorm{ y}_2$. To see this, suppose \cref{eqn-add-subtrct-event-case-2} holds and $\inlinenorm{y}_2 \leq R+\delta$. Then $L(y,x) \leq L_{R+\delta}$. If we now apply \cref{eqn-no-lipschitz-control} and this inequality,
\begin{equation}
0 \leq F(y)-\flb \leq F(x) - \flb - \dot F(x)^\intercal (y-x) + \frac{L_{R+\delta}}{2} \norm{y-x}_2^2,
\end{equation}
which contradicts \cref{eqn-add-subtrct-event-case-2}. Hence, in this case, $R+\delta < \inlinenorm{y}_2$. 

Using the triangle inequality, $R+\delta < \inlinenorm{y}_2 \leq \inlinenorm{x}_2 + \inlinenorm{y-x}_2 \leq R + \inlinenorm{y-x}_2$. That is, $1 \leq \inlinenorm{ y-x}_2/\delta \leq \inlinenorm{ y-x}_2^2/\delta^2$.

Hence,
\begin{align}
&\left[F(x) - \flb - \dot F(x)^\intercal (y-x) + \frac{L_{R+\delta}}{2} \norm{y-x}_2^2\right] [\pi_y(R)\pi_x(R) - \pi_x(R)] \nonumber\\
 &\quad \leq \delta\norm{\dot F(x)}_2 \frac{\norm{ y-x}_2}{\delta} \left[ \pi_x(R) - \pi_y(R)\pi_x(R)\right]\\
 &\quad \leq \delta \norm{\dot F(x)}_2 \frac{\norm{y-x}_2^2}{\delta^2}\left[ \pi_x(R) - \pi_y(R)\pi_x(R)\right]\\
 &\quad \leq \frac{G_R}{\delta}  \norm{ y-x}_2^2 \pi_x(R),
\end{align}
where in the last line we have used $\pi_x(R) - \pi_y(R) \pi_x(R) \leq \pi_x(R)$ as these are $\lbrace 0,1 \rbrace$-valued quantities.

Putting these two cases together in \cref{eqn-add-subtract-event}, we conclude,
\begin{equation}
\begin{aligned}
&[F(y) - \flb] \pi_y(R)\pi_x(R) \\
 &\quad \leq \left[F(x) - \flb - \dot F(x)^\intercal (y-x) + \left(\frac{L_{R+\delta}}{2} + \frac{G_R}{\delta} \right) \norm{y-x}_2^2\right] \pi_x(R).
\end{aligned}
\end{equation}
Letting $C_R = L_{R+\delta}/2 + G_R/\delta$, the conclusion follows.
\end{proof}
\begin{remark}
If we replace $(y,x)$ with $(x_{k+1},x_k)$ in the preceding result, it might be tempting to choose a $\delta$ that minimizes $C_R$ and then to use a standard approach to find a complexity result. However, this complexity result would only hold if all of the iterates remained within a radius $R$ of $0$, which, under \cref{assumption-lower-bound,assumption-gradient}, cannot be guaranteed apriori as shown by our construction in \cref{section-divergence-example}. Thus, a complexity result would only be appropriate if some additional information is known to guarantee a single Lipschitz constant (e.g., by knowing that the iterates remain bounded), in which case we would directly make use of \cref{equation-global-upper-bound-model} and we would have no use for \cref{result-recursion-gap}.
\end{remark}

We now apply \cref{result-recursion-gap} to the iterate sequence generated by gradient descent. To do so, we will make use of the following notation
\begin{equation}
\chi_{k}^0(R) = \begin{cases}
1 & \inlinenorm{x_{j}}_2 \leq R, ~j=0,\ldots,k \\
0 & \text{otherwise}.
\end{cases}
\end{equation}
That is, $\chi_k^0(R) = \pi_{x_0}(R)\pi_{x_1}(R)\cdots\pi_{x_k}(R)$. With this notation, we have the following simplification of \cref{result-recursion-gap} when applied to gradient descent.

\begin{corollary} \label{result-recursion-gap-simplified}
Suppose $F:\mathbb{R}^p \to \mathbb{R}$ satisfies \cref{assumption-lower-bound,assumption-gradient}. Let $x_0 \in \mathbb{R}^p$ and let $\lbrace x_k : k \in \mathbb{N} \rbrace$ be generated by \cref{eqn-recursion-iterate} satisfying \cref{property-spd,property-diminishing}. Then $\forall R \geq 0$, $\exists K \in \mathbb{N}$ such that for all $k \geq K$, 
\begin{equation}
[F(x_{k+1}) - \flb]\chi_{k+1}^0(R) \leq \left[F(x_k) - \flb - \frac{1}{2}\lambda_{\min}(M_k) \norm{\dot F(x_k)}_2^2 \right]\chi_{k}^0(R).
\end{equation}
\end{corollary}
\begin{proof}
By \cref{result-recursion-gap}, $\exists C_R > 0$ such that, for any $k+1 \in \mathbb{N}$,
\begin{equation}
\begin{aligned}
&[F(x_{k+1}) - \flb] \pi_{x_{k+1}}(R)\pi_{x_k}(R) \\
 &\quad \leq \left[F(x_k) - \flb - \dot F(x_k)^\intercal M_k \dot F(x_k) + C_R \norm{M_k \dot F(x_k)}_2^2\right] \pi_{x_k}(R),
\end{aligned}
\end{equation}
where we have made use of \cref{eqn-recursion-iterate} to replace $x_{k+1} - x_k$. If we now multiply both sides by the non-negative quantity $\pi_{x_0}(R)\cdots\pi_{x_{k-1}}(R)$, then
\begin{equation}
\begin{aligned}
&[F(x_{k+1}) - \flb] \chi_{k+1}^0(R)\\
 &\quad \leq \left[F(x_k) - \flb - \dot F(x_k)^\intercal M_k \dot F(x_k) + C_R \norm{M_k \dot F(x_k)}_2^2\right] \chi_{k}^0(R).
\end{aligned}
\end{equation}

The result follows if we show that $\exists K \in \mathbb{N}$ such that $\forall k \geq K$, 
\begin{equation} \label{eqn-recursion-gap-simplified-sufficient}
- \dot F(x_k)^\intercal M_k \dot F(x_k) + C_R \norm{ M_k \dot F(x_k)}_2^2 \leq -\frac{1}{2} \lambda_{\min}(M_k) \norm{ \dot F(x_k)}_2^2.
\end{equation}
To this end, we prove, if $M$ is symmetric positive definite with $\lambda_{\max}(M) < 1/(2C_R)$ then, for any $v \in \mathbb{R}^p$ with unit norm,
$- v^\intercal M v + C_R v^\intercal MM v \leq   - \frac{1}{2} \lambda_{\min}(M).$ Let $0 < \lambda_{\min}(M) = \lambda_p \leq \lambda_{p-1} \leq \cdots \leq \lambda_2 \leq \lambda_1 = \lambda_{\max}(M) < 1/(2C_R)$, where $\lambda_{\ell}$ denote the eigenvalues of $M$. Using the Schur Decomposition, there exists an orthogonal matrix $Q$ such that
$
- v^\intercal M v + C_R v^\intercal MM v = \sum_{\ell=1}^p ( - \lambda_\ell + C_R \lambda_{\ell}^2 ) w_\ell^2,
$
where $w_\ell$ is the $\ell^\mathrm{th}$ component of $Qv$ (note, $\inlinenorm{w}_2 = \inlinenorm{Qv}_2 = \inlinenorm{v}_2 = 1$). Since $\lambda_{\ell} < 1/(2C_R)$, it follows that $C_R \lambda_{\ell}^2 < \lambda_{\ell}/2$. Subtracting $\lambda_{\ell}$ from both sides, $-\lambda_{\ell} + C_R \lambda_{\ell}^2 < - \lambda_{\ell}/2 \leq -\lambda_{\min}(M)/2$. Thus,
\begin{equation}
- v^\intercal M v + C_R v^\intercal MM v \leq - \sum_{\ell=1}^p \frac{\lambda_\ell}{2} w_\ell^2 = -\frac{1}{2}v^\intercal M v \leq - \frac{\lambda_{\min}(M)}{2}.
\end{equation}
Since $\lambda_{\max}(M_k) \to 0$, there exists a $K \in \mathbb{N}$ such that $\forall k \geq K$, $\lambda_{\max}(M_k) \leq 1/(2C_R)$. Hence, there exists a $K$ such that $\forall k \geq K$, \cref{eqn-recursion-gap-simplified-sufficient} holds.
\end{proof}

\subsection{Applying Zoutendjik's Analysis Approach} \label{subsection-zoutendjik}

We now apply the recursive relationship established in \cref{result-recursion-gap-simplified} to study the objective and gradient using Zoutendjik's analysis method \cite{zoutendijk1960}. Recall, our main conclusion from the next result is that the limit suprmeum and limit infimum of the objective function evaluated at the iterates must tend to each other if the iterates persist in a region for long enough (even if they eventually escape), and the limit infimum of the gradient function evaluated at the iterates must tend to zero under similar circumstances. We stress that these conclusions are not the same as presupposing that the iterates remain in a bounded region.

\begin{theorem} \label{result-zoutendjik}
Suppose $F:\mathbb{R}^p \to \mathbb{R}$ satisfies \cref{assumption-lower-bound,assumption-gradient}. Let $x_0 \in \mathbb{R}^p$ and let $\lbrace x_k : k \in \mathbb{N} \rbrace$ be generated by \cref{eqn-recursion-iterate} satisfying \cref{property-spd,property-diverge,property-diminishing}. Then, for all $ R \geq 0$,
\begin{equation}
\lim_{k \to \infty} F(x_k)\chi_k^0(R) ~\text{exists and is finite, and}~
\liminf_{k \to \infty} \norm{ \dot F(x_k) }_2 \chi_k^0(R) = 0.
\end{equation}
If $\sup_{k} \inlinenorm{ x_k}_2 < \infty$ then $\lim_{k \to \infty} F(x_k)$ exists and is finite, and $\liminf_{k \to \infty} \inlinenorm{ \dot F(x_k) }_2 = 0$.
\end{theorem}
\begin{proof}
Let $R \geq 0$. The conditions of \cref{result-recursion-gap-simplified} are satisfied, and its conclusion is used freely herein. For the objective function, $\exists K \in \mathbb{N}$ such that for all $k \geq K$, $[ F(x_{k+1}) - \flb ] \chi_{k+1}^0(R) \leq [F(x_k) - \flb] \chi_k^0(R)$. As $\lbrace [ F(x_k) - \flb ] \chi_k^0(R) : k \geq K \rbrace$ is a nonincreasing sequence bounded from below, it converges. Now, if we further assume that $\sup_k \inlinenorm{x_k}_2 < \infty$, then there exists an $R > 0$ such that $\chi_k^0(R) = 1$ $\forall k+1 \in \mathbb{N}$. Hence, $\lim_{k \to \infty} F(x_k) - \flb$ exists and is finite.

For the gradient function, applying the conclusion of \cref{result-recursion-gap-simplified} and rearranging terms, for all $k \geq K$,
\begin{equation}
\frac{1}{2}\lambda_{\min}(M_k) \norm{ \dot F(x_k) }_2^2 \chi_k^0(R) \leq [F(x_{k}) - \flb] \chi_{k}^0(R) - [ F(x_{k+1}) - \flb ]\chi_{k+1}^0(R).
\end{equation}
Letting $j \geq K$ and using $F(x_{j+1}) - \flb \geq 0$, 
\begin{equation} \label{eqn-zoutendjik-telescoping-sum}
\sum_{k=K}^j \frac{1}{2}\lambda_{\min}(M_k) \norm{ \dot F(x_k) }_2^2 \chi_k^0(R) \leq [F(x_{K}) - \flb] \chi_{K}^0(R).
\end{equation}
Now for a contradiction, suppose $\exists c > 0$ such that $\liminf_{k \to \infty} \inlinenorm{ \dot F(x_k)}_2^2 \chi_k^0(R) > c$. Then, there exists a $K' > K$ such that 
\begin{equation}
\frac{c}{2} \sum_{k=K'}^j \lambda_{\min}(M_k) \leq \sum_{k=K'}^j \frac{1}{2}\lambda_{\min}(M_k) \norm{ \dot F(x_k) }_2^2 \chi_k^0(R) \leq [F(x_{K}) - \flb]\chi_{K}^0(R) < \infty.
\end{equation}
By \cref{property-diverge}, we have a contradiction. This part of the result follows for any $R \geq 0$.

Now, if $\sup_{k} \inlinenorm{x_k}_2 < \infty$ then there exists an $R > 0$ such that $\sup_{k} \inlinenorm{ x_k}_2 < R$. Therefore, $\chi_k^0(R) = 1$ for all $k+1 \in \mathbb{N}$, and, thus, the final part of the result follows.
\end{proof}

\begin{remark} \label{discussion-no-lipschitz-control}
Suppose we directly attempt to use Zoutendjik's analysis approach in \cref{eqn-no-lipschitz-control} with $y=x_{k+1}$ and $x=x_{k}$. We begin by rearranging \cref{eqn-no-lipschitz-control} and summing up to $j \in \mathbb{N}$ to conclude,
\begin{equation}
\sum_{k=0}^j \dot F(x_k)^\intercal M_k \dot F(x_k) - \frac{L(x_{k+1},x_k)}{2} \norm{M_k \dot F(x_k)}_2^2 \leq  F(x_0) - \flb.
\end{equation}
Thus, we conclude
\begin{equation}
\lim_{k \to \infty} \dot F(x_k)^\intercal M_k \dot F(x_k) - \frac{L(x_{k+1},x_k)}{2} \norm{M_k \dot F(x_k)}_2^2 = 0.
\end{equation}
Unfortunately, this conclusion does not imply that $\dot F(x_k) \to 0$ as $k \to \infty$. For instance, suppose as $k \to \infty$, $L(x_{k+1},x_k) \to \infty$. If $M_k = 2L(x_{k+1},x_k)^{-1} I$ for all $k$, then a straightforward substitution will show that the limit is satisfied yet $\dot F(x_k)$ does not have to be zero. Hence, using Zoutendjik's analysis method on this line of logic would not produce the desired conclusion. However, as shown in \cref{result-zoutendjik}, using Zoutendjik's analysis method on the conclusion of \cref{result-recursion-gap} is fruitful.
\end{remark}

\subsection{Convergence of the Gradient} \label{subsection-convergence-gradient}

One limitation of \cref{result-zoutendjik} is that it only provides for the limit infimum of the gradient function to be zero. Here, we will use \cref{property-conditioning} to conclude that the limit of the gradient function is zero.

\begin{theorem} \label{result-convergence-gradient}
Suppose $F:\mathbb{R}^p \to \mathbb{R}$ satisfies \cref{assumption-lower-bound,assumption-gradient}. Let $x_0 \in \mathbb{R}^p$ and let $\lbrace x_k : k \in \mathbb{N} \rbrace$ be generated by \cref{eqn-recursion-iterate} satisfying \cref{property-spd,property-diverge,property-diminishing,property-conditioning}. Then, for all $ R \geq 0$,
\begin{equation}
\lim_{k \to \infty} F(x_k)\chi_k^0(R) ~\text{exists and is finite, and}~
\lim_{k \to \infty} \norm{ \dot F(x_k) }_2 \chi_k^0(R) = 0.
\end{equation}
If $\sup_{k} \inlinenorm{ x_k}_2 < \infty$ then $\lim_{k \to \infty} F(x_k)$ exists and is finite, and $\lim_{k \to \infty} \inlinenorm{ \dot F(x_k) }_2 = 0$.
\end{theorem}
\begin{proof}
By \cref{result-zoutendjik}, we need only prove, for any $R \geq 0$, $\limsup_{k \to \infty} \inlinenorm{ \dot F(x_k) }_2\chi_{k}^0(R) = 0$. Fix $R \geq 0$. There are two cases.
\paragraph{Case 1} For some $K + 1 \in \mathbb{N}$, $\chi_{K}^0(R) = 0$. Then, $\chi_{k}^0(R) = 0$ for all $k \geq K$. The result follows.
\paragraph{Case 2} For all $k+1 \in \mathbb{N}$, $\chi_k^0(R) = 1$. In this case, $\inlinenorm{x_k}_2 \leq R$ for all $k+1 \in \mathbb{N}$. 
Let $L_R$ be the Lipschitz constant in the closed ball of radius $R$ around $0$ (see \cref{result-equivalent-lipschitz-definitions}), and let $G_R$ be the supremum of $\inlinenorm{ \dot F(x) }_2$ over all $x$ in the closed ball of radius $R$ around $0$. 

We now proceed in two steps. First, we show that for any $\epsilon > 0$, there exists a $K' \in \mathbb{N}$ such that $\forall k \geq K'$,
\begin{equation}
\left\vert \norm{ \dot F(x_{k+1}) }_2 - \norm{ \dot F(x_{k}) }_2 \right\vert < \frac{\epsilon}{4}.
\end{equation}
Then, we use a proof-by-contradiction to show that the $\limsup_{k \to \infty} \inlinenorm{ \dot F(x_k) }_2 \not > \epsilon$.

For the first part, let $\epsilon > 0$. Now,
\begin{align}
\left\vert \norm{ \dot F(x_{k+1}) }_2 - \norm{ \dot F(x_{k}) }_2 \right\vert 
&\leq \norm{ \dot F(x_{k+1}) - \dot F(x_k) }_2\\
&\leq L_R \norm{ x_{k+1} - x_k }_2\\
&\leq L_R \norm{ M_k \dot F(x_k) }_2\\
&\leq L_R G_R \lambda_{\max}(M_k).
\end{align}
By \cref{property-diminishing}, $\exists K' \in \mathbb{N}$ such that, $\forall k \geq K'$, $L_R G_R \lambda_{\max}(M_k) < \epsilon / 4$. 

Suppose now $\limsup_{k \to \infty} \inlinenorm{ \dot F(x_k) }_2 > \epsilon$. Let $u_0 = \min \lbrace k > \max \lbrace K, K' \rbrace : \inlinenorm{ \dot F(x_k)}_2 > \epsilon \rbrace$, where $K$ is given by \cref{result-recursion-gap-simplified}. By \cref{result-zoutendjik}, we can now define the following three subsequences of $\mathbb{N}$ for all $i \in \mathbb{N}$:
\begin{remunerate}
\item $j_i = \min \lbrace t > u_{i-1} : \inlinenorm{ \dot F(x_t) }_2 < \epsilon / 2 \rbrace$.
\item $u_i = \min \lbrace t > j_i : \inlinenorm{ \dot F(x_t) }_2 > \epsilon \rbrace$.
\item $\ell_i = \min \lbrace t \in [j_i,u_i) : \inlinenorm{ \dot F(x_s) }_2 > \epsilon /2, s=t+1,\ldots,u_i \rbrace$.
\end{remunerate}
Importantly, since $\ell_i > K'$ for all $i \in \mathbb{N}$ and $\inlinenorm{ \dot F(x_{\ell_i + 1})}_2 > \epsilon/2$, $\inlinenorm{ \dot F(x_{\ell_i})}_2 > \epsilon / 4$. Using these subsequences and the same steps from the first part,
\begin{align}
\frac{\epsilon}{2} 
&< \norm{ \dot F(x_{u_i}) }_2 - \norm{ \dot F(x_{\ell_i})}_2\\
&< \sum_{t = \ell_i}^{u_i-1} \norm{ \dot F(x_{t+1})}_2 - \norm{ \dot F(x_t)}_2\\
&< \sum_{t = \ell_i}^{u_i-1} L_R \lambda_{\max}(M_k) \norm{ \dot F(x_k) }_2\\
&< \frac{\epsilon}{4} \sum_{t=\ell_i}^{u_i-1} L_R \lambda_{\max}(M_k) \left( \frac{4 \norm{\dot F(x_k)}_2}{\epsilon} \right)\\
&< \frac{\epsilon}{4} \sum_{t=\ell_i}^{u_i-1} L_R \lambda_{\max}(M_k) \left( \frac{4 \norm{\dot F(x_k)}_2}{\epsilon} \right)^2,
\end{align}
where we made use of $\epsilon/4 < \inlinenorm{ \dot F(x_s) }_2$ for $s=\ell_i,\ldots,u_i-1$ in the ultimate line. Simplifying and applying \cref{property-conditioning}, $\forall i \in \mathbb{N}$,
\begin{equation} \label{equation-gradient-gap-contradiction}
\frac{\epsilon^2}{8L_R \kappa} < \sum_{ t=\ell_i}^{u_i - 1} \lambda_{\min}(M_k) \norm{ \dot F(x_k)}_2^2. 
\end{equation}
Summing both sides over $i \in \mathbb{N}$, the left hand side diverges while the right hand side is bounded by \cref{eqn-zoutendjik-telescoping-sum}. Hence, we have a contradiction and the conclusion follows.
\end{proof}

From this proof, we might question whether it is necessary to use \cref{property-conditioning} in order to replace $\lambda_{\max}(M_k)$ with $\lambda_{\min}(M_k)$ in \cref{equation-gradient-gap-contradiction}. We provide a concrete example where our reasoning faces difficulty if \cref{property-conditioning} is not used. As the example below shows, it is possible to relax \cref{property-conditioning} if the sequence $\lbrace M_k \rbrace$ eventually has common invariant subspaces, but we do not pursue this here.

\begin{example} \label{example-condition-control}
Let $F: \mathbb{R}^2 \to \mathbb{R}$ be
\begin{equation}
F(x) = \frac{1}{2} (\cmp{x}{1})^2 + \frac{1}{10} (\cmp{x}{2})^2,
\end{equation}
where $\cmp{x}{i}$ is the $i^\mathrm{th}$ component of $x$. Consider now $x_0$ such that $\cmp{x_0}{1} = 0$ and $\cmp{x_0}{2} = 1$. In order to violate \cref{property-conditioning}, let 
\begin{equation}
M_k = \frac{1}{5}\begin{bmatrix}
(k+1)^{-1/2} & 0 \\
0 & (k+1)^{-1}
\end{bmatrix}, ~k + 1 \in \mathbb{N}.
\end{equation}
When we apply gradient descent, $\cmp{x_k}{1} = 0$ for all $k \in \mathbb{N}$, and $\cmp{x_k}{2} > 0.8(k+1)^{-1/5}$ \cite[p.1578]{nemirovski2009robust}. Then, $\inlinenorm{ \dot F(x_k) }_2 > 0.16 (k+1)^{-1/5}$. Now, we have,
$
\lambda_{\max}(M_k) \inlinenorm{ \dot F(x_k)}_2^2 > 0.01 (k+1)^{-9/10},
$
which produces a divergent series, whereas $\lambda_{\min}(M_k)$ in place of $\lambda_{\max}(M_k)$ would produce a convergent series.\footnote{To show convergence, we need an upper bound on the rate of convergence of $\cmp{x_k}{2}$, which is on the order of $(k+1)^{-1/5}$ \cite[see][p.1578]{nemirovski2009robust}.}
\end{example}

\subsection{Topological Properties of the Iterates} \label{subsection-topological-iterates}
We now turn our attention to the asymptotic behavior of the iterates in the bounded regime. We will make use of the closure of subsequential limits (see \cref{result-closure-of-subsequential-limits}) and a fact about the density of subsequential limits of a decaying sequence (see \cref{result-density-of-subsequential-limits}). We state the main result in \cref{result-asymptotics-of-iterates}.

\begin{theorem} \label{result-asymptotics-of-iterates}
Suppose $F:\mathbb{R}^p \to \mathbb{R}$ satisfies \cref{assumption-lower-bound,assumption-gradient}. Let $x_0 \in \mathbb{R}^p$ and let $\lbrace x_k : k \in \mathbb{N} \rbrace$ be generated by \cref{eqn-recursion-iterate} satisfying \cref{property-spd,property-diverge,property-diminishing,property-conditioning}. If $\sup_{k} \inlinenorm{x_k}_2 < \infty$ and we let $\mathcal{C}$ denote the subsequential limits of $\lbrace x_k : k+1 \in \mathbb{N} \rbrace$, then
\begin{remunerate}
\item $\mathcal{C}$ is closed;
\item $\forall z \in \mathcal{C}$, $\dot F(z) = 0$;
\item $\mathcal{C}$ is connected;
\item $\mathcal{C}$ does not contain an open set; and
\item Either $|\mathcal{C}| = 1$ or $|\mathcal{C}| = \infty$.
\end{remunerate}
\end{theorem}
\begin{proof}
The first statement follows from \cref{result-closure-of-subsequential-limits}. For the second statement, if $z \in \mathcal{C}$ then there is a subsequence $\lbrace x_{k_j} : j \in \mathbb{N} \rbrace$ such that $\lim_{j} x_{k_j} = z$. By the continuity of $x \mapsto \dot F(x)$ (see \cref{assumption-gradient}), $\dot F(z) = \lim_{j} \dot F(x_{k_j})$. The limit on the right hand side is zero by \cref{result-convergence-gradient}.

For the third statement, recall that $\mathcal{C}$ is bounded by hypothesis and $\mathcal{C}$ is closed by the first statement. Hence, $\mathcal{C}$ is compact. Suppose that $\mathcal{C}$ is not connected. Then, there are two disjoint open sets, $O_1$ and $O_2$, whose union contains $\mathcal{C}$ and whose individual intersections with $\mathcal{C}$ are non-empty. We denote the intersections of $O_1$ and $O_2$ with $\mathcal{C}$ by $\mathcal{C}_1$ and $\mathcal{C}_2$, respectively. We now proceed in three steps. First, we verify that $\mathcal{C}_1$ and $\mathcal{C}_2$ are closed, and, consequently, compact. Second, we use compactness to show that the distance between $\mathcal{C}_1$ and $\mathcal{C}_2$ is strictly larger than zero. Third, we use the diminishing step sizes and \cref{result-density-of-subsequential-limits} to derive a contradiction.

Suppose $\mathcal{C}_1$ is not closed. Let $z$ be a limit point of $\mathcal{C}_1$ that is not in $\mathcal{C}_1$. Then $z \in \mathcal{C}$, which implies that $z \in \mathcal{C}_2 \subset O_2$. There is a sequence of points in $\mathcal{C}_1$ contained in an arbitrarily small neighborhood of $z$, which implies $\mathcal{C}_1 \cap O_2 \neq \emptyset$, which is a contradiction. Hence, $\mathcal{C}_1$ is closed. The same argument shows $\mathcal{C}_2$ is closed.

Since $\mathcal{C}_1$ and $\mathcal{C}_2$ are closed and bounded, they are compact. Now, $(z_1,z_2) \mapsto \inlinenorm{ z_1 - z_2 }_2$ is a continuous function. Hence, this function applied to $\mathcal{C}_1 \times \mathcal{C}_2$ must achieve its minimum at some points $z_1^* \in \mathcal{C}_1$ and $z_2^* \in \mathcal{C}_2$. If $z_1^* = z_2^*$, then $O_1 \cap O_2 \neq \emptyset$ which is a contradiction. Hence, $z_1^* \neq z_2^*$ so the distance between any points in $\mathcal{C}_1$ and $\mathcal{C}_2$ is at least $\inlinenorm{z_1^* - z_2^*}_2 > 0$.

Define a function $g: \mathbb{R}^p \to \mathbb{R}_{\geq 0 }$ such that $g(x) = \inf_{w \in \mathcal{C}_1} \inlinenorm{ x - w}_2$. Then, $g(z) = 0$ for any $z \in \mathcal{C}_1$ and $g(z) \geq \inlinenorm{z_1^* - z_2^*}_2$ for $z \in \mathcal{C}_2$. Hence, $\liminf_k g(x_k) = 0$ and $\limsup_k g(x_k) \geq \inlinenorm{z_1^* - z_2^*}$. We now verify, $\lim_k g(x_{k+1}) - g(x_k) = 0$, and apply \cref{result-density-of-subsequential-limits} to derive a contradiction. For any $k \in \mathbb{N}$, there exists a $w_k \in \mathcal{C}_1$ such that $g(x_k) = \inlinenorm{x_k - w_k}_2$. Hence,
\begin{align}
g(x_{k+1}) - g(x_k) &= \inf_{w \in \mathcal{C}_1} \inlinenorm{ x_{k+1} - w}_2 - \inlinenorm{ x_k - w_k}_2 \\ 
				   &\leq \inlinenorm{ x_{k+1} - w_k}_2 - \inlinenorm{ x_k - w_k}_2 \\
				   &\leq \inlinenorm{ x_{k+1} - x_k}_2\\
				   &\leq \inlinenorm{ M_k \dot F(x_k)}_2.
\end{align}
Note, $\inlinenorm{ M_k \dot F(x_k)}_2 \leq \lambda_{\max}(M_k) G_R$, where $G_R = \sup_{x : \inlinenorm{x}_2 \leq R} \inlinenorm{ \dot F(x)}_2$ and $R = \sup_{k} \inlinenorm{x_k}_2 < \infty$. Since $\lambda_{\max}(M_k) \to 0$, then $g(x_{k+1}) - g(x_k) \to 0$. Hence, by \cref{result-density-of-subsequential-limits}, there is a subsequence, $\lbrace g(x_{k_j}) : j \in \mathbb{N} \rbrace$ that converges to, say, $\inlinenorm{z_1^* - z_2^*}_2 / 2$. Consequently, $\lbrace x_{k_j}: j \in \mathbb{N} \rbrace$ is a bounded sequence, and it has a subsequence that converges to a point $z^*$ such that $\inf_{w \in \mathcal{C}_1} \inlinenorm{ z^* - w}_2 = \inlinenorm{z_1^* - z_2^*}_2/2$. Hence, $z^*$ is a subsequential limit but it is not in either $\mathcal{C}_1$ nor $\mathcal{C}_2$, which is a contradiction. Thus, $\mathcal{C}$ is connected.

For the fourth statement, suppose $\mathcal{C}$ contains an open set $O$. Let $z \in O$. Since $z$ is a limit point of a subsequence, there exists an $k \in \mathbb{N}$ such that $x_k \in O$. By the first statement, $\dot F(x_k) = 0$, which, by \cref{eqn-recursion-iterate}, implies $x_{j} = x_k$ for all $j \geq k$. Hence, $\mathcal{C}$ is the singleton, $\lbrace x_k \rbrace$, which is a contradiction. Thus, $\mathcal{C}$ cannot contain an open set.

For the final statement, recall that $\mathcal{C}$ is connected. This implies that $\mathcal{C}$ cannot contain a finite number of points other than a single point. So, either $|\mathcal{C}|=1$ or $|\mathcal{C}| = \infty$.
\end{proof}

\section{The Divergence Regime} \label{section-divergence-example}
\Cref{result-zoutendjik} leaves open the possibility that the iterates can diverge. Of course, this divergence regime is possible even under the stricter assumption of global Lipschitz continuity of the gradient. Importantly, under global Lipschitz continuity of the gradient, when the iterates diverge, the objective function still converges to a finite quantity and the gradient function converges to zero \cite[Proposition 1.2.4]{bertsekas2016}. For example, the globally Lipschitz continuous function, $F(x) = \exp(-x^2)$, achieves its minimum as the iterates diverge; and, in this divergence regime, the objective function converges to zero and the gradient function converges to zero.

Under our more realistic assumption of local Lipschitz continuity, will the objective function converge to a finite quantity and will the gradient function converge to zero when the iterates diverge? Unfortunately, the answer is no. In this section, we will construct several examples that show the extreme behaviors which can occur in the divergence regime when only local Lipschitz continuity is assumed. Of note, we construct an example in which catastrophic divergence can occur: the iterates diverge, the objective function diverges to infinity, and the gradient norm remains uniformly bounded away from zero. We show this construction here. Our remaining constructions are found as specified in the following \cref{table-divergence-examples}. We underscore that our constructions can be used to generate objective functions on which gradient descent can have other interesting behaviors that we do not explicitly construct here (e.g., the limit infimum and limit supremum of the gradient function being distinct).

\begin{table}[H]
\centering
\caption{Summary of Counterexamples for the Divergence Regime.} \label{table-divergence-examples}
\begin{tabular}{lp{4.5in}} \toprule
\textbf{Reference} & \textbf{Summary} \\ \midrule
This section & A case for which the iterates of gradient descent will produce objective function values that diverge and gradient function values that are uniformly bounded away from zero. \\
\Cref{section-divergence-objective-limit-nonexistence} & A case for which the iterates of gradient descent will produce objective function values whose limit supremum is infinity and whose limit infimum is zero, while the gradient function values remain bounded away from zero. \\
\Cref{subsection-divergence-gradient-zero} & A case for which the iterates of gradient descent will produce objective function values that diverge and the gradient function tends to zero. \\ \bottomrule
\end{tabular}
\end{table}

\subsection{Construction of the Objective Function}
Let $\lbrace m_k : k+1 \in \mathbb{N}\rbrace$ be a sequence of scalars such that $m_k > 0$, $\sum_k m_k = \infty$, and $m_k \to 0$ as $k \to \infty$. Define $S_0 = 0$ and $S_{k+1} = \sum_{j=0}^k m_k$ for all integers $k \geq 0$. 

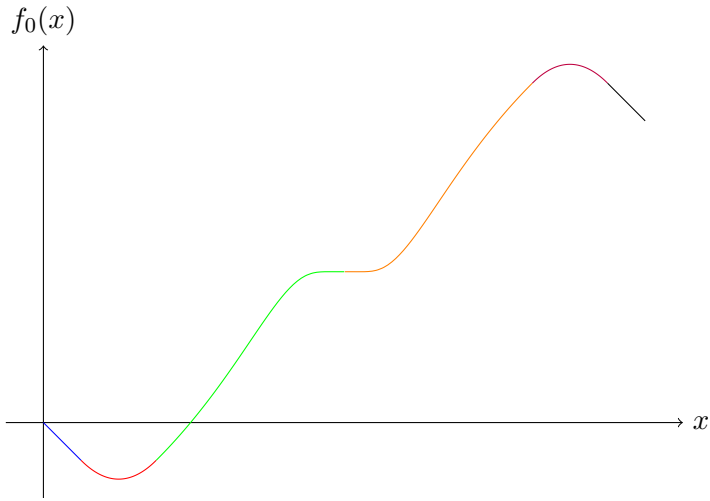
\begin{figure}[h]
\centering
\begin{tikzpicture}[scale=1]
	\def\m{8.0}
	\draw[->] (-0.5, 0) -- (\m+0.5, 0) node[right] {$x$};
	\draw[->] (0,-1) -- (0, \m/2 + 1.0) node[above] {$f_0(x)$};
	\draw[scale=1, domain=0:\m/16, smooth, variable=\x, blue] plot ({\x}, {-\x});
	\draw[scale=1, domain=\m/16:3*\m/16, smooth, variable=\x, red]  plot ({\x}, {(8/\m)*(\x - \m/8)^2 - 3*\m/32});
	\draw[scale=1, domain=3*\m/16:\m/2, smooth, variable=\x, green] plot ({\x}, { (-5*\m/16)*exp( (5*\m/16)/( \x - \m/2) + 1) + \m/4 });
	\draw[scale=1, domain=\m/2 + \m/1000: 13*\m/16, smooth, variable=\x, orange] plot ({\x}, { (5*\m/16)*exp( -(5*\m/16)/( \x - \m/2) + 1) + \m/4 });
	\draw[scale=1, domain=13*\m/16:15*\m/16, smooth, variable=\x, purple]  plot ({\x}, {-(8/\m)*(\x - 7*\m/8)^2 + 19*\m/32});
	\draw[scale=1, domain=15*\m/16:\m, smooth, variable=\x, black] plot ({\x}, {-\x + 3*\m/2});
\end{tikzpicture}
\caption{Plot of $f_0(x)$ with $m_0 = 8.0$ with each component shown in a different color.}
\label{plot-divergent-example-component}
\end{figure}

For the objective function, define $F: \mathbb{R} \to \mathbb{R}$ by
\begin{equation} \label{eqn-divergence-objective}
F(x) = \begin{cases}
-x & x \leq 0 \\
f_j(x) & x \in (S_j, S_{j+1}], ~\forall j+1 \in \mathbb{N},
\end{cases}
\end{equation}
where $\lbrace f_j: (S_j, S_{j+1}] \to \mathbb{R} : j + 1\in\mathbb{N} \rbrace$ are defined iteratively as follows. Let
\begin{equation} \label{eqn-divergence-component-base}
f_0(x) = \begin{cases}
-x & x \in (0, \frac{m_0}{16}) \\
\frac{8}{m_0}(x- \frac{m_0}{8})^2 - \frac{3m_0}{32} & x \in [\frac{m_0}{16}, \frac{3m_0}{16}) \\
- \frac{5m_0}{16} \exp\left( \frac{5m_0/16}{x - m_0/2 } + 1 \right) + \frac{m_0}{4} & x \in [\frac{3m_0}{16}, \frac{ m_0}{2}) \\
\frac{m_0}{4} & x = \frac{m_0}{2} \\
\frac{5m_0}{16} \exp \left( -\frac{5m_0/16}{x - m_0/2} + 1 \right) + \frac{m_0}{4} & x \in (\frac{m_0}{2}, \frac{13m_0}{16}) \\
\frac{-8}{m_0}(x - \frac{7m_0}{8})^2 + \frac{19m_0}{32} & x \in [\frac{13 m_0}{16}, \frac{15m_0}{16}) \\
-x + \frac{3 m_0}{2} & x \in [\frac{15 m_0}{16}, m_0],
\end{cases}
\end{equation}
which is plotted for a particular choice of $m_0$ in \cref{plot-divergent-example-component}. Now, for $j \in \mathbb{N}$, let $x' = x - S_j$ and let
\begin{equation} \label{eqn-divergence-component-general}
f_j(x) = \begin{cases}
-x' + f_{j-1}(S_j) & x' \in (0, \frac{m_j}{16}) \\
\frac{8}{m_j}(x' - \frac{m_j}{8})^2 - \frac{3m_j}{32} + f_{j-1}(S_j) & x' \in [\frac{m_j}{16}, \frac{3m_j}{16}) \\
- \frac{5m_j}{16} \exp\left( \frac{5m_j/16}{x' - m_j/2} + 1\right) + \frac{m_j}{4} + f_{j-1}(S_j) & x' \in [\frac{3m_j}{16},\frac{m_j}{2}) \\
\frac{m_j}{4} + f_{j-1}(S_j) & x' = \frac{m_j}{2} \\
\frac{5m_j}{16} \exp\left( \frac{-5m_j/16}{x' - m_j/2} + 1\right) + \frac{m_j}{4} + f_{j-1}(S_j) & x' \in (\frac{m_j}{2},\frac{13m_j}{16}) \\
\frac{-8}{m_j}(x' - \frac{7m_j}{8})^2 + \frac{19m_j}{32} + f_{j-1}(S_j) & x' \in [\frac{13m_j}{16}, \frac{15m_j}{16}) \\
-x' + \frac{3m_j}{2} + f_{j-1}(S_j) & x' \in [\frac{15m_j}{16}, m_j].
\end{cases}
\end{equation}

\subsection{Properties of the Objective Function}

We show that $F: \mathbb{R} \to \mathbb{R}$ as defined in \cref{eqn-divergence-objective}--\cref{eqn-divergence-component-general} satisfies \cref{assumption-lower-bound,assumption-gradient}. We begin by proving that each component, $f_j: (S_j, S_{j+1}] \to \mathbb{R}$, satisfies \cref{assumption-lower-bound,assumption-gradient} on its domain. 

\begin{remark} \label{remark-continuous-extension}
Below, we define the continuous extension of $f_j$ on $[S_j, S_{j+1}]$ by the value of $f_j(x)$ on $(S_j, S_{j+1}]$ and by $\lim_{x \downarrow S_j} f_j(x)$ for the point $x= S_j$. Moreover, at the ends of the interval, we use differentiability and the corresponding notation, $\dot f_j$, to mean the one-sided derivatives.
\end{remark}
\begin{proposition} \label{result-divergence-component-properties}
The continuous extension $f_0: (S_0, S_1] \to \mathbb{R}$ $($as defined in \cref{eqn-divergence-component-base}$)$ to $[S_0,S_1]$ is continuous on its domain, bounded from below by $-3m_0/32$, differentiable on its domain with $\dot f_0(S_0) = \dot f_0(S_1) = -1$, and its derivative is locally Lipschitz continuous. Similarly, the continuous extension of $f_j: (S_j, S_{j+1}] \to \mathbb{R}$ $($as defined in \cref{eqn-divergence-component-general}$)$ to $[S_j,S_{j+1}]$ is continuous on its domain, bounded from below by $f_{j-1}(S_j) - 3m_j/32$, differentiable on its domain with $\dot f_j(S_j) = -1$ and $\dot f_j(S_{j+1}) = -1$, and its derivative is locally Lipschitz continuous.
\end{proposition}
\begin{proof}
We only look at an arbitrary $j \in \mathbb{N}$ as the proof is identical for $f_0$. Moreover, since $f_j$ is equal to its reflection across the vertical axis $x = S_j + m_j/2$ followed by a reflection over the horizontal axis $y = m_j/4 + f_{j-1}(S_j)$, it is enough to show continuity and differentiability on $[S_j, S_j + m_j/2]$. 

To establish continuity, we need to show that the left sided limits of $f_j$ agree with the function value at $S_j + \delta m_j/16$ for $\delta = 1, 3, 8$. Starting with $\delta = 1$, $\lim_{x \uparrow S_j + m_j/16} -x + S_j + f_{j-1}(S_j) = -m_j/16 + f_{j-1}(S_j)$. By direct substitution,
\begin{equation}
f_j(S_j + m_j/16) = \frac{8}{m_j}\left( \frac{m_j}{16} \right)^2 - \frac{3m_j}{32} + f_{j-1}(S_j) = \frac{m_j}{32} - \frac{3m_j}{32} + f_{j-1}(S_j).
\end{equation}
Hence, the left limit agrees with the function value at $\delta = 1$. For $\delta = 3$, the symmetry and continuity of the quadratic function implies $\lim_{ x \uparrow S_j + 3m_j/16} f_j(x) = -m_j/16 + f_{j-1}(S_j)$. By direct substitution,
\begin{equation}
f_j(S_j + 3m_j/16) = -\frac{5m_j}{16} \exp\left( \frac{5m_j/16}{-5m_j/16} + 1 \right) + \frac{4m_j}{16} + f_{j-1}(S_j) = -\frac{m_j}{16} + f_{j-1}(S_j).
\end{equation}
Hence, the left limit agrees with the function value at $\delta = 3$. For $\delta = 8$,
\begin{equation}
\lim_{ x \uparrow S_j + m_j/2} -\frac{5m_j}{16} \exp\left( \frac{5m_j/16}{x - S_j - m_j/2} + 1 \right) + \frac{m_j}{4} + f_{j-1}(S_j) = \frac{m_j}{4} + f_{j-1}(S_j), 
\end{equation} 
which is just $f_j(S_j + m_j/2)$. Hence, the continuous extension of $f_j$ is continuous on $[S_j, S_{j+1}]$. 

We now compute the derivatives of the components of $f_j$ on $[S_j, S_{j+1}]$ with the convention of assigning the one-sided derivative to the component function that includes its end point. Let $x' = x - S_j$. 
\begin{equation} \label{eqn-divergence-component-derivatives}
\dot f_j(x) =\begin{cases}
-1 & x' \in [0, \frac{m_j}{16}) \\
\frac{16}{m_j}(x' - \frac{m_j}{8}) & x' \in [\frac{m_j}{16},\frac{3m_j}{16}) \\
\left( \frac{5m_j/16}{x' - m_j/2} \right)^2 \exp\left( \frac{5m_j/16}{x' - m_j/2} + 1 \right) & x' \in [\frac{3m_j}{16}, \frac{m_j}{2}) \\
\left( \frac{5m_j/16}{x' - m_j/2} \right)^2 \exp\left( \frac{-5m_j/16}{x' - m_j/2} + 1 \right) & x' \in (\frac{m_j}{2}, \frac{13m_j}{16}) \\
-\frac{16}{m_j}(x' - \frac{7m_j}{8}) & x' \in [\frac{13m_j}{16}, \frac{15m_j}{16}) \\
-1 & x' \in [\frac{15m_j}{16}, m_j].
\end{cases}
\end{equation}
In order to extend these component derivatives to the continuous extension of $f_j$, we need to verify that the left hand limits of the derivatives agree with the right side derivatives at $S_j + \delta m_j/16$ for $\delta=1,3$ of \cref{eqn-divergence-component-derivatives}, and we need to verify that the left hand limit of the derivative is $0$ at $\delta = 8$ of \cref{eqn-divergence-component-derivatives}. Starting with $\delta = 1$, the left hand limit is $-1$ and, by direct calculation,
\begin{equation}
\frac{16}{m_j}\left(S_j + \frac{m_j}{16} - S_j - \frac{m_j}{8} \right) = -1.
\end{equation}
For $\delta = 3$, the left hand limit is
\begin{equation}
\lim_{ x \uparrow S_{j} + {3m_j/16} } \frac{16}{m_j}\left(x - S_j - \frac{m_j}{8} \right) = 1,
\end{equation}
and a direct evaluation of the third component \cref{eqn-divergence-component-derivatives} is 
\begin{equation}
\left( \frac{5m_j/16}{S_{j} + {3m_j/16} - S_j - m_j/2} \right)^2 \exp\left( \frac{5m_j/16}{S_{j} + {3m_j/16} - S_j - m_j/2} + 1 \right) = 1.
\end{equation}
For $\delta = 8$, we need to check that the left hand limit is zero, which can be confirmed by checking that the argument of the exponential term goes to $-\infty$ as $x \uparrow S_j + 8m_j/16$.

Overall, the derivative of the continuous extension of $f_j$ is well defined at every point on its interval, is continuous, and is given by (with $x' = x - S_j$)
\begin{equation} \label{eqn-divergence-component-derivative-definition}
\dot f_j(x) = \begin{cases}
-1 & x' \in [0, \frac{m_j}{16}) \\
\frac{16}{m_j}(x' - \frac{m_j}{8}) & x' \in [\frac{m_j}{16}, \frac{3m_j}{16}) \\
\left( \frac{5m_j/16}{x' - m_j/2} \right)^2 \exp\left( \frac{5m_j/16}{x' - m_j/2} + 1 \right) & x' \in [\frac{3m_j}{16}, \frac{m_j}{2}) \\
0 & x' = m_j/2 \\
\left( \frac{5m_j/16}{x' - m_j/2} \right)^2 \exp\left( \frac{-5m_j/16}{x' - m_j/2} + 1 \right) & x' \in (\frac{m_j}{2}, \frac{13m_j}{16}) \\
-\frac{16}{m_j}(x' - \frac{7m_j}{8}) & x' \in [\frac{13m_j}{16}, \frac{15m_j}{16}) \\
-1 & x' \in [\frac{15m_j}{16}, m_j].
\end{cases}
\end{equation}

Using this derivative, we can calculate the lower bound for the function. By the derivative of the extension of $f_j$, \cref{eqn-divergence-component-derivative-definition}, we see that the function is decreasing only on $[S_j, S_j + m_j/8]$ and $[S_j + 7m_j/8, S_{j+1}]$. Moreover, $f_j(S_j + m_j/8) = -3m_j/32 + f_{j-1}(S_j)$ and $f_j(S_{j+1}) = m_j/2 + f_{j-1}(S_j)$. Thus, the lower bound of the extension of $f_j$ is as stated.

Our last step is to verify the local Lipschitz continuity of the derivative. It is easy to verify that within its interval, the components are twice continuously differentiable. As a result, we can use \cref{result-hessian-implies-lipschitz}. Similarly, if we define the second derivative at $S_j + m_j/2$ to be $0$, we can verify that the objective is twice continuously differentiable at $S_j + m_j/2$. Then, we can use \cref{result-hessian-implies-lipschitz} again. To conclude, we need to examine what happens around the points $S_j + \delta m_j/16$ for $\delta = 1, 3$.

Starting with $\delta = 1$, consider the points $S_j + m_j/16 - \epsilon_1 m_j/16$ and $S_j + m_j/16 + \epsilon_2 m_j/16$ for $\epsilon_1, \epsilon_2 > 0$ sufficiently small. Then the difference in the derivatives at these points divided by the distance between the points is
\begin{align}
\frac{\left\vert \frac{16}{m_j}\left(-\frac{m_j}{16} + \epsilon_2 \frac{m_j}{16}\right) + 1  \right\vert }{(\epsilon_2 + \epsilon_1) m_j/16} = \frac{\epsilon_2}{(\epsilon_2 + \epsilon_1) m_j/16} \leq \frac{16}{m_j}.
\end{align} 
Therefore, we conclude that the derivative is locally Lipschitz near $S_j + m_j/16$. For $\delta=3$, we compute the same ratio at the points $S_j + 3m_j/16 - \epsilon_1 m_j/16$ and $S_j + 3m_j/16 + 5 \epsilon_2 m_j/16$ for $\epsilon_1, \epsilon_2 \in [0,1/4]$ where at most either $\epsilon_1$ or $\epsilon_2$ is zero. The ratio of the difference in the derivatives and the points is
\begin{align}
\frac{\left\vert \frac{1}{(1 - \epsilon_2)^2} \exp\left( \frac{\epsilon_2}{\epsilon_2 - 1}\right) -1 + \epsilon_1 \right\vert}{(5 \epsilon_2 + \epsilon_1) m_j/16} \leq \frac{\epsilon_2 + \epsilon_2^2/2 + \epsilon_1}{(5 \epsilon_2 + \epsilon_1) m_j/16} \leq \frac{16}{m_j}.
\end{align}
Therefore, we conclude that the derivative is locally Lipschitz near $S_j + 3m_j/16$. 
\end{proof}

With this calculation complete, we can now verify that $F$ satisfies \cref{assumption-lower-bound,assumption-gradient}.

\begin{proposition} \label{result-divergence-objective-properties}
The function $F : \mathbb{R} \to \mathbb{R}$ is continuous and differentiable on its domain; the function $F$ is lower bounded; the derivative of the function $F$ is locally Lipschitz continuous; $F(S_j) = S_j/2$; and $\dot F(S_j) = -1$ for all $j + 1 \in \mathbb{N}$. And, the derivative of the function $F$ is not globally Lipschitz continuous.
\end{proposition}
\begin{proof}
Using \cref{result-divergence-component-properties}, we can verify continuity of $F$ by checking the continuity of $F$ at $x = S_j$ for all $j$. Note, $F(S_j) = f_{j-1}(S_j)$. We need to verify that the right side limit converges to $f_{j-1}(S_j)$. That is,
\begin{equation}
\lim_{ x \downarrow S_j } F(x) = \lim_{x \downarrow S_j} f_j(x) = \lim_{x \downarrow S_j} -x + S_j + f_{j-1}(S_j) = f_{j-1}(S_j). 
\end{equation} 
Thus, $F$ is continuous at all $S_j$ for $j \in \mathbb{N}$. We check $x = S_0 = 0$ as well. $F(0) = 0$ by definition. Moreover,
\begin{equation}
\lim_{ x\downarrow 0 } F(x) = \lim_{x \downarrow 0} f_0(x) = \lim_{ x \downarrow 0} -x = 0.
\end{equation}
Hence, $F$ is continuous on its domain.

For differentiability, we similarly need only verify the differentiability of $F$ at $x=S_j$ for all $j$. By  \cref{result-divergence-component-properties}, it follows that the derivative at each $S_j$ is $-1$ from the left and the right for $j \in \mathbb{N}$. Hence, the derivative exists at each $S_j$ for $j \in \mathbb{N}$. Moreover, the left hand derivative at $S_0$ is $-1$, which agrees with the right hand derivative. Hence, $F$ is continuously differentiable on its domain. Since the derivative is constant in a small neighborhood of $S_j$, it is Lipschitz continuous in this region. Finally, $\dot F(S_j) = -1$ for all $j+1 \in \mathbb{N}$. 

To show that $F$ is lower bounded, we will first calculate the values of $F(S_j)$. We proceed by induction. $F(S_0) = 0 = S_0/2$. Now suppose the statement holds up to $j$. Then $F(S_{j}) = S_{j}/2$. By construction, $f_{j-1}(S_j) = F(S_{j}) = S_j/2$. Now, 
\begin{align}
F(S_{j+1}) = f_j(S_{j+1}) &= -S_{j+1} + S_j + \frac{3m_j}{2} + f_{j-1}(S_j) \\
&= -S_j - m_j + S_j + \frac{3m_j}{2} + \frac{S_j}{2} = \frac{S_{j+1}}{2}.
\end{align}
To show the lower bound property, recall that $f_j(x) \geq f_{j-1}(S_j) - 3m_j/32 = F(S_{j}) - 3m_j/32 = S_j/2 - 3m_j/32$. Since $S_j \to \infty$ and $m_j \to 0$ by construction, $F(x) \geq \inf_{j} S_j/2 - 3m_j/32 > -\infty$.

Finally, we verify that $F$ is not globally Lipschitz continuous. For a contradiction, suppose there exists an $L > 0$ such that for any $x,x' \in \mathbb{R}$, $|\dot F(x) - \dot F(x') | \leq L |x - x'|$. Since $m_j \to 0$, there exists $j \in \mathbb{N}$ such that $L m_j / 2 < 1$. Then,
$| \dot F(S_j + m_j/2) - \dot F(S_j) | \leq Lm_j/2 < 1$. However, by \cref{result-divergence-component-properties}, $\dot F(S_j + m_j/2) = \dot f_j(S_j + m_j/2) = 0$ and $\dot F(S_j) = \dot f_j(S_j) = -1$, which implies $| \dot F(S_j + m_j/2) - \dot F(S_j) | = 1$, which is a contradiction.
\end{proof}

\subsection{Properties of Gradient Descent on the Objective Function}

We are now ready to show that gradient descent with diminishing step sizes generates iterates that diverge, and whose objective function diverges and  gradient function remains bounded away from zero.

\begin{proposition} \label{result-divergence}
Let $\lbrace m_k : k+1 \in \mathbb{N} \rbrace$ be any positive sequence such that $\sum_{k} m_k$ diverges and $m_k \to 0$. Define $F : \mathbb{R} \to \mathbb{R}$ as in \cref{eqn-divergence-objective}. Suppose $x_0 = 0$ and let $\lbrace x_k : k \in \mathbb{N} \rbrace$ be generated according to \cref{eqn-recursion-iterate} with $M_k = m_k I$ for all $k+1 \in \mathbb{N}$. Then, $\lbrace M_k \rbrace$ satisfies \cref{property-spd,property-diverge,property-diminishing}. Moreover, (a) $\lim_k x_k = \infty$, (b) $\lim_k F(x_k) = \infty$, and (c) $\lim_k |\dot F(x_k)| = 1$. 
\end{proposition}
\begin{proof}
To prove the result, we need only show that $x_k = S_k$ where we recall that $S_0 = 0$ and $S_k = \sum_{j=0}^{k-1} m_j$. For $k=0$, $x_0 = 0 = S_0$. Suppose this holds up to $k$. Then, by \cref{result-divergence-objective-properties},
\begin{equation}
x_{k+1} = x_k - M_k \dot F(x_k) = S_k - m_k (-1) = S_{k+1}.
\end{equation}
Now, since $S_k$ diverges, the iterates diverge (part $(a)$). Moreover, by \cref{result-divergence-objective-properties}, since $F(x_k) = F(S_k) = S_k/2$, the objective function also diverges (part $(b)$). Finally, by \cref{result-divergence-objective-properties}, $\dot F(x_k) = \dot F(S_k) = -1$ (part $(c)$).
\end{proof}

In summary, as the example from \cref{result-divergence} and the example of $F(x) = \exp(-x^2)$ show, under our assumptions about the objective function and properties of gradient descent, we cannot conclude anything additional about the objective behavior of the function or the  gradient in the regime where the iterates generated by gradient descent with diminishing step sizes diverge.

\section{Conclusion} \label{section-conclusion}
In this paper, we have analyzed the global behavior of gradient descent with diminishing step sizes for differentiable nonconvex functions whose gradients are only locally Lipschitz continuous. To the best of our knowledge, we have provided the most general convergence analysis of gradient descent with diminishing step sizes. Specifically, we have shown that the iterates cannot produce erratic behavior in the objective function nor gradient function when they persist in a region for sufficiently long, even if they eventually escape. We also construct specific examples to show the types of erratic behaviors which can occur when the iterates escape off to infinity. Our analysis has also raised a number of interesting questions with varying degrees of practical interest.
\begin{remunerate}
\item Is there a notion of continuity on the gradients that is appropriate for data science yet more restrictive than \cref{assumption-gradient} for which \cref{result-zoutendjik} or \cref{result-convergence-gradient} hold uniformly over the family of functions specified by this notion of continuity?
\item Is there a choice of step sizes that ensures the subsequential limit points of the iterates is a set that is a singleton?
\item Is there a function class that is necessary and sufficient to avoid the divergence regime and the corresponding erratic behaviors for gradient descent with diminishing step size?
\end{remunerate}

%\section*{Acknowledgments}
%The authors thank the detailed and insightful feedback and criticism from the associate editors and reviewers,
%which has substantially improved the quality of this work.

\bibliographystyle{siamplain}
\bibliography{gd}

\begin{thebibliography}{10}

\bibitem{armijo1966}
{\sc L.~Armijo}, {\em Minimization of functions having lipschitz continuous
  first partial derivatives}, Pacific Journal of mathematics, 16 (1966),
  pp.~1--3.

\bibitem{beck2017first}
{\sc A.~Beck}, {\em First-order methods in optimization}, SIAM, 2017.

\bibitem{benaim1999}
{\sc M.~Bena{\"\i}m}, {\em Dynamics of stochastic approximation algorithms}, in
  Seminaire de probabilites XXXIII, Springer, 1999, pp.~1--68.

\bibitem{benveniste2012}
{\sc A.~Benveniste, M.~M{\'e}tivier, and P.~Priouret}, {\em Adaptive algorithms
  and stochastic approximations}, vol.~22, Springer Science \& Business Media,
  2012.

\bibitem{bertsekas2016}
{\sc D.~Bertsekas}, {\em Nonlinear Programming}, vol.~4, Athena Scientific,
  2016.

\bibitem{bittanti1990convergence}
{\sc S.~Bittanti, P.~Bolzern, and M.~Campi}, {\em Convergence and exponential
  convergence of identification algorithms with directional forgetting factor},
  Automatica, 26 (1990), pp.~929--932.

\bibitem{bottou2010}
{\sc L.~Bottou}, {\em Large-scale machine learning with stochastic gradient
  descent}, in Proceedings of COMPSTAT'2010, Springer, 2010, pp.~177--186.

\bibitem{cao2003exponential}
{\sc L.~Cao and H.~M.~Schwartz}, {\em Exponential convergence of the kalman
  filter based parameter estimation algorithm}, International Journal of
  Adaptive Control and Signal Processing, 17 (2003), pp.~763--783.

\bibitem{cartis2011adaptiveI}
{\sc C.~Cartis, N.~I. Gould, and P.~L. Toint}, {\em Adaptive cubic
  regularisation methods for unconstrained optimization. part i: motivation,
  convergence and numerical results}, Mathematical Programming, 127 (2011),
  pp.~245--295.

\bibitem{cartis2011adaptiveII}
{\sc C.~Cartis, N.~I. Gould, and P.~L. Toint}, {\em Adaptive cubic
  regularisation methods for unconstrained optimization. part ii: worst-case
  function-and derivative-evaluation complexity}, Mathematical programming, 130
  (2011), pp.~295--319.

\bibitem{cauchy1847}
{\sc A.~Cauchy et~al.}, {\em M{\'e}thode g{\'e}n{\'e}rale pour la
  r{\'e}solution des systemes d’{\'e}quations simultan{\'e}es}, Comp. Rend.
  Sci. Paris, 25 (1847), pp.~536--538.

\bibitem{curry1944}
{\sc H.~B. Curry}, {\em The method of steepest descent for non-linear
  minimization problems}, Quarterly of Applied Mathematics, 2 (1944),
  pp.~258--261.

\bibitem{curtis2020adaptive}
{\sc F.~E. Curtis and K.~Scheinberg}, {\em Adaptive stochastic optimization: A
  framework for analyzing stochastic optimization algorithms}, IEEE Signal
  Processing Magazine, 37 (2020), pp.~32--42.

\bibitem{curtis2019stochastic}
{\sc F.~E. Curtis, K.~Scheinberg, and R.~Shi}, {\em A stochastic trust region
  algorithm based on careful step normalization}, Informs Journal on
  Optimization, 1 (2019), pp.~200--220.

\bibitem{defossez2020simple}
{\sc A.~D{\'e}fossez, L.~Bottou, F.~Bach, and N.~Usunier}, {\em A simple
  convergence proof of adam and adagrad}, arXiv preprint arXiv:2003.02395,
  (2020).

\bibitem{dommel1968}
{\sc H.~W. Dommel and W.~F. Tinney}, {\em Optimal power flow solutions}, IEEE
  Transactions on power apparatus and systems,  (1968), pp.~1866--1876.

\bibitem{du2019}
{\sc S.~Du, J.~Lee, H.~Li, L.~Wang, and X.~Zhai}, {\em Gradient descent finds
  global minima of deep neural networks}, in International conference on
  machine learning, PMLR, 2019, pp.~1675--1685.

\bibitem{du2017}
{\sc S.~S. Du, C.~Jin, J.~D. Lee, M.~I. Jordan, A.~Singh, and B.~Poczos}, {\em
  Gradient descent can take exponential time to escape saddle points}, Advances
  in neural information processing systems, 30 (2017).

\bibitem{du2018}
{\sc S.~S. Du, X.~Zhai, B.~Poczos, and A.~Singh}, {\em Gradient descent
  provably optimizes over-parameterized neural networks}, arXiv preprint
  arXiv:1810.02054,  (2018).

\bibitem{fort1996}
{\sc J.-C. Fort and G.~Pages}, {\em Convergence of stochastic algorithms: From
  the kushner--clark theorem to the lyapounov functional method}, Advances in
  applied probability, 28 (1996), pp.~1072--1094.

\bibitem{grapiglia2022adaptive}
{\sc G.~N. Grapiglia and G.~F. Stella}, {\em An adaptive trust-region method
  without function evaluations}, Computational Optimization and Applications,
  82 (2022), pp.~31--60.

\bibitem{gratton2022convergence}
{\sc S.~Gratton, S.~Jerad, and P.~L. Toint}, {\em Convergence properties of an
  objective-function-free optimization regularization algorithm, including an
  $o(\epsilon^{3/2})$ complexity bound}, arXiv preprint arXiv:2203.09947,
  (2022).

\bibitem{gratton2022first}
{\sc S.~Gratton, S.~Jerad, and P.~L. Toint}, {\em First-order
  objective-function-free optimization algorithms and their complexity}, arXiv
  preprint arXiv:2203.01757,  (2022).

\bibitem{gratton2022parametric}
{\sc S.~Gratton, S.~Jerad, and P.~L. Toint}, {\em Parametric complexity
  analysis for a class of first-order adagrad-like algorithms}, arXiv preprint
  arXiv:2203.01647,  (2022).

\bibitem{hadamard1908}
{\sc J.~Hadamard}, {\em M{\'e}moire sur le probl{\`e}me d'analyse relatif {\`a}
  l'{\'e}quilibre des plaques {\'e}lastiques encastr{\'e}es}, vol.~33,
  Imprimerie nationale, 1908.

\bibitem{iyengar2013}
{\sc G.~Iyengar and A.~K.~C. Ma}, {\em Fast gradient descent method for
  mean-cvar optimization}, Annals of Operations Research, 205 (2013),
  pp.~203--212.

\bibitem{jin2018}
{\sc C.~Jin, P.~Netrapalli, and M.~I. Jordan}, {\em Accelerated gradient
  descent escapes saddle points faster than gradient descent}, in Conference On
  Learning Theory, PMLR, 2018, pp.~1042--1085.

\bibitem{johnstone1982exponential}
{\sc R.~M. Johnstone, C.~R. Johnson~Jr, R.~R. Bitmead, and B.~D. Anderson},
  {\em Exponential convergence of recursive least squares with exponential
  forgetting factor}, Systems \& Control Letters, 2 (1982), pp.~77--82.

\bibitem{karimi2016}
{\sc H.~Karimi, J.~Nutini, and M.~Schmidt}, {\em Linear convergence of gradient
  and proximal-gradient methods under the polyak-{\l}ojasiewicz condition}, in
  Joint European Conference on Machine Learning and Knowledge Discovery in
  Databases, Springer, 2016, pp.~795--811.

\bibitem{ke1998class}
{\sc X.~Ke and J.~Han}, {\em A class of nonmonotone trust region algorithms for
  unconstrained optimization problems}, Science in China Series A: Mathematics,
  41 (1998), pp.~927--932.

\bibitem{lee2019}
{\sc J.~Lee, L.~Xiao, S.~Schoenholz, Y.~Bahri, R.~Novak, J.~Sohl-Dickstein, and
  J.~Pennington}, {\em Wide neural networks of any depth evolve as linear
  models under gradient descent}, Advances in neural information processing
  systems, 32 (2019).

\bibitem{lee2016}
{\sc J.~D. Lee, M.~Simchowitz, M.~I. Jordan, and B.~Recht}, {\em Gradient
  descent only converges to minimizers}, in Conference on learning theory,
  PMLR, 2016, pp.~1246--1257.

\bibitem{lemarechal2012}
{\sc C.~Lemar{\'e}chal}, {\em Cauchy and the gradient method}, Doc Math Extra,
  251 (2012), p.~10.

\bibitem{ljung1977}
{\sc L.~Ljung}, {\em Analysis of recursive stochastic algorithms}, IEEE
  transactions on automatic control, 22 (1977), pp.~551--575.

\bibitem{mertikopoulos2020}
{\sc P.~Mertikopoulos, N.~Hallak, A.~Kavis, and V.~Cevher}, {\em On the almost
  sure convergence of stochastic gradient descent in non-convex problems},
  arXiv preprint arXiv:2006.11144,  (2020).

\bibitem{more1983recent}
{\sc J.~J. Mor{\'e}}, {\em Recent developments in algorithms and software for
  trust region methods}, Springer, 1983.

\bibitem{nemirovski2009robust}
{\sc A.~Nemirovski, A.~Juditsky, G.~Lan, and A.~Shapiro}, {\em Robust
  stochastic approximation approach to stochastic programming}, SIAM Journal on
  optimization, 19 (2009), pp.~1574--1609.

\bibitem{nocedal2006}
{\sc J.~Nocedal and S.~Wright}, {\em Numerical optimization}, Springer Science
  \& Business Media, 2006.

\bibitem{oymak2019}
{\sc S.~Oymak and M.~Soltanolkotabi}, {\em Overparameterized nonlinear
  learning: Gradient descent takes the shortest path?}, in International
  Conference on Machine Learning, PMLR, 2019, pp.~4951--4960.

\bibitem{parkum1992recursive}
{\sc J.~Parkum, N.~K. Poulsen, and J.~Holst}, {\em Recursive forgetting
  algorithms}, International Journal of Control, 55 (1992), pp.~109--128.

\bibitem{patel2020}
{\sc V.~Patel}, {\em Stopping criteria for, and strong convergence of,
  stochastic gradient descent on bottou-curtis-nocedal functions}, Mathematical
  Programming,  (2021), pp.~1--42.

\bibitem{patel2021global}
{\sc V.~Patel, B.~Tian, and S.~Zhang}, {\em Global convergence and stability of
  stochastic gradient descent}, arXiv preprint arXiv:2110.01663,  (2021).

\bibitem{patel2021}
{\sc V.~Patel and S.~Zhang}, {\em Stochastic gradient descent on nonconvex
  functions with general noise models}, arXiv preprint arXiv:2104.00423,
  (2021).

\bibitem{reddi2016}
{\sc S.~Reddi, S.~Sra, B.~Poczos, and A.~J. Smola}, {\em Proximal stochastic
  methods for nonsmooth nonconvex finite-sum optimization}, Advances in neural
  information processing systems, 29 (2016), pp.~1145--1153.

\bibitem{reddi2016stochastic}
{\sc S.~J. Reddi, A.~Hefny, S.~Sra, B.~Poczos, and A.~Smola}, {\em Stochastic
  variance reduction for nonconvex optimization}, in International conference
  on machine learning, PMLR, 2016, pp.~314--323.

\bibitem{ward2020adagrad}
{\sc R.~Ward, X.~Wu, and L.~Bottou}, {\em Adagrad stepsizes: Sharp convergence
  over nonconvex landscapes}, The Journal of Machine Learning Research, 21
  (2020), pp.~9047--9076.

\bibitem{wu2018wngrad}
{\sc X.~Wu, R.~Ward, and L.~Bottou}, {\em Wngrad: Learn the learning rate in
  gradient descent}, arXiv preprint arXiv:1803.02865,  (2018).

\bibitem{zhang2006superlinearly}
{\sc J.~Zhang, Y.~Wang, and X.~Zhang}, {\em Superlinearly convergent
  trust-region method without the assumption of positive-definite hessian},
  Journal of optimization theory and applications, 129 (2006), pp.~201--218.

\bibitem{zhang2007trust}
{\sc J.~Zhang, L.~Wu, and X.~Zhang}, {\em A trust region method for
  optimization problem with singular solutions}, Applied Mathematics and
  Optimization, 56 (2007), pp.~379--394.

\bibitem{zoutendijk1960}
{\sc G.~Zoutendijk}, {\em Methods of feasible directions: a study in linear and
  non-linear programming}, Elsevier, 1960.

\end{thebibliography}

\appendix
\section{Equivalent Definitions for Local Lipschitz Continuity}
\begin{lemma} \label{result-equivalent-lipschitz-definitions}
A function $G: \mathbb{R}^p \to \mathbb{R}^p$ is locally Lipschitz continuous if and only if for every compact set $\mathcal{C} \subset \mathbb{R}^p$, there exists an $L \geq 0$ such that
\begin{equation}
\frac{ \norm{ G(y) - G(z)}_2 }{\norm{y - z}_2 } \leq L, ~\forall y,z\in \mathcal{C}.
\end{equation}
\end{lemma}
\begin{proof}
Suppose $G$ is locally Lipschitz continuous. Suppose for a contradiction, there exists a compact set for which no such $L$ exists. Then for every $\ell \in \mathbb{N}$, we can find a pair $y_\ell, z_\ell\in \mathcal{C}$ such that
\begin{equation}
\frac{ \norm{ G(y_\ell) - G(z_\ell)}_2 }{\norm{ y_\ell - z_\ell}_2} > \ell.
\end{equation}
By compactness, there exists a subsequence $\lbrace \ell_k : k \in \mathbb{N} \rbrace$ and $y,z\in \mathcal{C}$ such that $y_{\ell_k} \to y$ and $z_{\ell_k} \to z$ as $k \to \infty$. If $\norm{y - z}_2 > 0$, then, for $k \in \mathbb{N}$ sufficiently large,
\begin{equation}
\frac{ \norm{ G(y_{\ell_k}) - G(z_{\ell_k})}_2 }{\norm{ y_\ell - z_\ell}_2} \leq \frac{2 \sup_{x \in \mathcal{C}} \norm{G(x)}_2 }{0.5 \norm{ y - z}_2} < \infty,
\end{equation}
which is a contradiction. Hence, $\norm{y - z}_2 = 0$; that is, $y = z$. This also provides a contradiction as $G$ is locally Lipschitz continuous at $y = z$ and so for $k \in \mathbb{N}$ sufficiently large, $y_{\ell_k}$ and $z_{\ell_k}$ would be inside of $\mathcal{N}$ from \cref{definition-local-lipschitz}.

For the other direction of the result: for any point $x \in \mathbb{R}^p$ and any open ball containing $x$, we can take the closure of this open ball to generate a compact set $\mathcal{C}$. The result follows.
\end{proof}

\section{Continuous Hessians Implies Local Lipschitz Continuity}
\begin{lemma} \label{result-hessian-implies-lipschitz}
Suppose $F$ is twice continuously differentiable for all $x \in \mathbb{R}^p$. Then $\dot F(x)$ is locally Lipschitz continuous.
\end{lemma}
\begin{proof}
Let $\ddot F(x)$ denote the Hessian of $F$. Then, by assumption, $\inlinenorm{ \ddot F(x)}_2$ is a continuous function and it is bounded over any compact region. By Taylor's theorem, for any $x,y \in \mathbb{R}^p$,
$\dot F(x) - \dot F(y) = \int_{0}^1 \ddot{F}(y + t (x - y) ) (x -y) dt$.
Let $K \subset \mathbb{R}^p$ be compact. By continuity and compactness, there exists an $L$ for $K$ such that $\inlinenorm{ \ddot F(x)}_2 \leq L$ for all $x \in K$. Hence, by H\"{o}lder's inequality, for any $x, y \in K$,
$\inlinenorm{ \dot F(x) - \dot F(y) } \leq L\inlinenorm{x - y}_2$.
As $K$ is arbitrary, the result follows.
\end{proof}

\section{Some Properties of Subsequential Limits}
\begin{lemma} \label{result-closure-of-subsequential-limits}
Let $\lbrace a_n : n \in \mathbb{N} \rbrace \subset \mathbb{R}^p$. Let $\mathcal{C}$ be the set of its subsequential limits. Then $\mathcal{C}$ is closed.
\end{lemma}
\begin{proof}
Let $z$ be a limit point of $\mathcal{C}$. Then, we can construct a sequence $\lbrace z_k : k \in \mathbb{N} \rbrace \subset \mathcal{C}$ such that for every $K \in \mathbb{N}$ and for all $k \geq K$, $\inlinenorm{ z_k - z }_2 \leq 2^{-K-1}$. Moreover, since $z_k \in \mathcal{C}$, $\exists n_k \in \mathbb{N}$ such that $\inlinenorm{ a_{n_k} - z_k}_2 \leq 2^{-k-1}$. Let $\epsilon > 0$ and let $K \in \mathbb{N}$ such that $2^{-K} < \epsilon$. Then, $\forall k \geq K$, $\inlinenorm{ a_{n_k} - z}_2 \leq \inlinenorm{ a_{n_k} - z_k}_2 + \inlinenorm{ z_k - z}_2 \leq 2^{-K} < \epsilon$. Hence, $z = \lim_{k} a_{n_k} \in \mathcal{C}$. 
\end{proof}
\begin{lemma} \label{result-density-of-subsequential-limits}
Let $\lbrace a_n : n \in \mathbb{N} \rbrace \subset \mathbb{R}$ such that $\liminf_{n} a_n$ and $\limsup_n a_n$ are finite. If $\lim_{n} a_{n+1} - a_n = 0$, then for any $z \in [\liminf_{n} a_n, \limsup_{n} a_n]$, there is a subsequence of $\lbrace a_n : n \in \mathbb{N} \rbrace$ that converges to $z$.
\end{lemma}
\begin{proof}
We begin by showing that any closed interval strictly between the limit infimum and limit supremum contains a subsequential limit. Let $r_1 < r_2$ such that $\liminf_{n} a_n < r_1$ and $r_2 < \limsup_{n} a_n$.
If there exists an infinite subsequence $\lbrace a_{n_k} : k \in \mathbb{N} \rbrace \subset [r_1,r_2]$, then sequential compactness implies that $\lbrace a_{n_k} : k \in \mathbb{N} \rbrace$ has a subsequence which converges in $[r_1,r_2]$. Suppose now, $\exists K \in \mathbb{N}$ such that $\forall n \geq K$, $a_n \not\in [r_1,r_2]$. Since the $\liminf_{n} a_n < r_1 < r_2 < \limsup_{n} a_n$, there exists a subsequence $\lbrace a_{n_k}: k \in \mathbb{N} \rbrace$ such that $a_{n_k} < r_1$ and $r_2 < a_{n_k+1}$. However, this is a contradiction since $a_{n_k + 1} - a_{n_k} \to 0$ as $k \to \infty$. Hence, there is always a subsequence in any closed interval between $\liminf_{n} a_n$ and $\limsup_n a_n$.

We have that if $z$ is either the limit infimum or limit supremum then there is a subsequence of $\lbrace a_n : n \in \mathbb{N} \rbrace$ that converges to this value. So take $\liminf_{n} a_n < z < \limsup_{n} a_n$. We now proceed by induction. Let $z_0 = \liminf_{n} a_n$. There is a subsequence that converges to a point in $[0.5(z + z_0), z]$. Let $z_1$ be this limit. If $z\neq z_1$, then $|z_1 - z| \leq 2^{-1} (z - z_0)$ and we define $z_2$ as the subsequential limit in $[0.5(z + z_1, z]$. If $z=z_1$ then we stop. Suppose we proceed by induction such that $\lbrace z_j : j=1,\ldots,k \rbrace$ are subsequential limits such that $|z_j - z| \leq 2^{-j}(z - z_0)$. If $z\neq z_k$, then we can find $z_{k+1}$ as the limit of a subsequence in $[0.5(z + z_k), z]$, which we denote $z_{k+1}$. Moreover, $|z_{k+1} - z| \leq 2^{-k-1}(z-z_0)$. If we never terminate at $z$ for some $k \in \mathbb{N}$, then $\lbrace z_k : k \in \mathbb{N} \rbrace$ is a sequence of subsequential limits converging to $z$. By \cref{result-closure-of-subsequential-limits}, $z$ is a subsequential limit.
\end{proof}

\section{Divergence Regime: Nonexistence of Objective Function Limit} \label{section-divergence-objective-limit-nonexistence}
Here, we use as similar construction for \cref{result-divergence} to construct an objective function $F$ such that when gradient descent is applied to this objective function with a specific initialization, $\limsup_k F(x_k) = \infty$, $\liminf_k F(x_k) = 0$ and $|\dot F(x_k) | = 1$ for all $k$. We proceed in three general steps corresponding to each subsection below.

\subsection{Objective Function Target Values}

Let $\lbrace m_k : k+1 \in \mathbb{N}\rbrace$ be a sequence of scalars such that $m_k > 0$, $\sum_k m_k = \infty$, and $m_k \to 0$ as $k \to \infty$. Define $S_0 = 0$ and $S_{k+1} = \sum_{j=0}^k m_k$ for all integers $k \geq 0$. We will now construct a sequence $\lbrace O_k : k+1 \in \mathbb{N} \rbrace$ which will serve as target values for each iterate of our objective function.

\begin{remunerate}
\item Let $O_0 = 0$. For convenience, let $u_0 = \ell_0 = 0$.
\item Let $\ell_1 = 1 + \min \lbrace k \geq 0 : O_0 + \frac{1}{2} \sum_{j=0}^k m_j > 1 \rbrace$. From the divergence of $\sum_k m_k$, it is clear that such an $\ell_1$ is finite. Define $O_{k} = O_0 + \frac{1}{2}\sum_{j=0}^{k - 1 } m_j$ for $k \in [1,\ell_1] \cap \mathbb{N}$.
\item Let $u_1 = 1 + \min \lbrace k \geq \ell_1 : O_{\ell_1} - \sum_{j=\ell_1}^k m_j < 0 \rbrace$. Again, from the divergence of $\sum_k m_k$, $u_1$ is finite. Define $O_{k} = O_{\ell_1} - \sum_{j=\ell_1}^{k - 1} m_j$ for $k \in [\ell_1+1,u_1]\cap \mathbb{N}$.
\item For $t \in \mathbb{N}$, let $\ell_{t+1} = 1 + \min \lbrace k \geq u_t : O_{u_t} + \frac{1}{2} \sum_{j=u_t}^k m_j > t+1 \rbrace$. From the divergence of $\sum_{k} m_k$, $\ell_{t+1}$ is finite if $u_t$ is finite. Define $O_{k} = O_{u_t} + \frac{1}{2} \sum_{j=u_t}^{k-1} m_j$ for $k \in [u_t+1,\ell_{t+1}] \cap \mathbb{N}$.
\item For $t \in \mathbb{N}$, let $u_{t+1} = 1 + \min \lbrace k \geq \ell_{t+1} : O_{\ell_{t+1}} - \sum_{j=\ell_{t+1}}^k m_j < 0 \rbrace$. From the divergence of $\sum_{k} m_k$, $u_{t+1}$ is finite if $\ell_{t+1}$ is finite. Define $O_{k} = O_{\ell_{t+1}} - \sum_{j = \ell_{t+1}}^{k-1} m_j$ for $k \in [\ell_{t+1} + 1,\ldots,u_{t+1}] \cap \mathbb{N}$.
\end{remunerate}

We point out several facts about the sequence $\lbrace O_k : k+1 \in \mathbb{N} \rbrace$. First, $\lim_t O_{\ell_{t}} = \infty$ by construction. Second, we verify, $\lim_{t} O_{u_t} = 0$. By construction, $O_{u_t-1} > 0$ and $0 > O_{u_t} = O_{u_t-1} - m_{u_t-1} \geq - m_{u_t-1}$. Since $m_{u_t-1} \to 0$ as $t \to \infty$, $\liminf_t O_{u_t} = 0$. In turn, the limit of the sequence exists and is zero. Third, we verify, $\limsup_k O_k = \infty$ and $\liminf_k O_k = 0$. For any $k \in \mathbb{N}_{>\ell_1}$, there exists a $t \in \mathbb{N}$ such that $k \in [\ell_t+1,u_t]$ or $k \in [u_t+1,\ell_{t+1}]$. If $k \in [\ell_t+1,u_t]$, then $O_k \in [O_{u_t}, O_{\ell_t}]$. If $k \in [u_{t}+1,\ell_{t+1}]$, then $O_k \in [O_{u_t}, O_{\ell_{t+1}}]$. Hence, the third fact holds because of the first two.

\subsection{Construction of the Objective Function}

With these sequences established, we now state our objective function.
\begin{equation} \label{eqn-divergence-objective-example}
F(x) = \begin{cases}
-x & x\leq 0, \\
\tilde f_0(x) & x \in (0, S_{\ell_1}], \\
\tilde f_t(x) & x \in (S_{\ell_t}, S_{\ell_{t+1}}], ~\forall t \in \mathbb{N},
\end{cases}
\end{equation}
where
\begin{equation}
\tilde f_0(x) = \begin{cases}
f_j(x) & x \in (S_j, S_{j+1}], ~ j \in \lbrace 0,\ldots, \ell_1-1 \rbrace;
\end{cases}
\end{equation}
\begin{equation} \label{eqn-divergence-component-objective-example}
\tilde f_t(x) = \begin{cases}
O_{\ell_t}- (x - S_{\ell_t}) & x \in (S_{\ell_t}, S_{u_t}] \\
f_j(x) & x \in (S_j, S_{j+1}], ~ j \in \lbrace u_t,\ldots, \ell_{t+1}-1 \rbrace;
\end{cases}
\end{equation}
and
\begin{equation} 
f_j(x) = \begin{cases}
-x' + O_j & x' \in (0, \frac{m_j}{16}) \\
\frac{8}{m_j}(x' - \frac{m_j}{8})^2 - \frac{3m_j}{32} + O_j & x' \in [\frac{m_j}{16}, \frac{3m_j}{16}) \\
- \frac{5m_j}{16} \exp\left( \frac{5m_j/16}{x' - m_j/2} + 1\right) + \frac{m_j}{4} + O_{j} & x' \in [\frac{3m_j}{16},\frac{m_j}{2}) \\
\frac{m_j}{4} + O_j & x' = \frac{m_j}{2} \\
\frac{5m_j}{16} \exp\left( \frac{-5m_j/16}{x' - m_j/2} + 1\right) + \frac{m_j}{4} + O_j & x' \in (\frac{m_j}{2},\frac{13m_j}{16}) \\
\frac{-8}{m_j}(x' - \frac{7m_j}{8})^2 + \frac{19m_j}{32} + O_j & x' \in [\frac{13m_j}{16}, \frac{15m_j}{16}) \\
-x' + \frac{3m_j}{2} + O_j & x' \in [\frac{15m_j}{16}, m_j].
\end{cases}
\end{equation}
with $x' = x - S_j$.

\subsection{Properties of the Objective Function}

Here, we verify, \cref{eqn-divergence-objective-example} satisfies \cref{assumption-lower-bound,assumption-gradient}.
We need to verify certain properties of $\tilde f_t(x)$, which we do now.

\begin{proposition} \label{result-divergence-component-objective-example}
Let $t+1 \in \mathbb{N}$. The continuous extension of $\tilde f_t : (S_{\ell_t}, S_{\ell_{t+1}}] \to \mathbb{R}$, \cref{eqn-divergence-component-objective-example}, to $[S_{\ell_t}, S_{\ell_{t+1}}]$ is
\begin{remunerate}
\item continuous on $[S_{\ell_t},S_{\ell_{t+1}}]$ with values $O_{\ell_t}$, $O_{u_t}$ and $O_{\ell_{t+1}}$ at points $S_{\ell_t}$, $S_{u_t}$ and $S_{\ell_{t+1}}$, respectively;
\item bounded from below by $\min \lbrace O_j - \frac{3m_j}{32} : j = u_t,\ldots, s_{\ell_{t+1} - 1} \rbrace$;
\item differentiable on $[S_{\ell_t}, S_{\ell_{t+1}}]$ with the one-sided derivatives being $-1$ at the end points of the interval;
\item locally Lipschitz continuous.
\end{remunerate}
\end{proposition}
\begin{proof}
We note that \cref{eqn-divergence-component-objective-example} has several components. The $f_j(x)$ are the same as those defined by \cref{eqn-divergence-component-base} but shifted vertically by a constant. Hence, by \cref{result-divergence-component-properties}, the continuous extension of $f_j(x)$ to $[S_{j},S_{j+1}]$ is continuous; bounded from below by $O_j - \frac{3m_j}{32}$; differentiable with the one-sided derivatives being $-1$ on the end points of the interval; and locally Lipschitz continuous.

We use these facts to show the remaining properties of $\tilde f_t(x)$. First, to verify continuity, we need only verify that the components agree at the points $x \in \lbrace S_{u_t}, S_{u_t+1},\ldots,S_{\ell_{t+1}-1} \rbrace$. When $x = S_{u_t}$, 
\begin{align}
\tilde f_t(S_{u_t}) = O_{\ell_t} + (S_{u_t} - S_{\ell_t}) = O_{\ell_t} + \left( \sum_{k=0}^{u_t-1} m_k - \sum_{k=0}^{\ell_t-1} m_k \right) = O_{\ell_t} + \sum_{k=\ell_t}^{u_t-1} = O_{u_t}.
\end{align}
Moreover,
\begin{align}
\lim_{x \downarrow S_{u_t}} \tilde f_t(x) = \lim_{x \downarrow S_{u_t}} f_{u_t}(x) = \lim_{x \downarrow S_{u_t}} -(x - S_{u_t}) + O_{u_t} = O_{u_t}.
\end{align}
Hence, the evaluation of $\tilde f_t(x)$ at $S_{u_t}$ agrees with its limit from the right. For the remaining points, let $j \in \lbrace u_t + 1,\ldots, \ell_{t+1} - 1 \rbrace$. Then,
\begin{equation}
\tilde f_t(S_j) = f_{j-1}(S_j) = -(S_j - S_{j-1}) + \frac{3m_{j-1}}{2} + O_{j-1} = \frac{m_{j-1}}{2} + O_{u_t} + \frac{1}{2}\sum_{k=u_t}^{j-1} m_k = O_j.
\end{equation}
Moreover,
\begin{equation}
\lim_{x \downarrow S_j} \tilde f_t(S_j) = \lim_{x \downarrow S_j} f_j(S_j) = \lim_{x \downarrow S_j} -(x - S_j) + O_j = O_j.
\end{equation}
Hence, $\tilde f(x)$ is continuous. Moreover, we have also shown that the continuous extension of $\tilde f(x)$ has the stated values at $x \in \lbrace S_{u_t}, S_{u_t+1},\ldots,S_{\ell_{t+1}-1} \rbrace$.

For the lower bound, we have that $\tilde f(x) \geq O_{u_t}$ for $x \in (S_{\ell_t}, S_{u_t}]$. By \cref{result-divergence-component-properties}, each $f_j(x) \geq O_j - \frac{3m_j}{32}$. Hence, the lower bound follows.

We now verify differentiability. By the properties of a linear function and \cref{result-divergence-component-properties}, each component of $\tilde f_t(x)$ is differentiable on its domain. We must check that these derivatives agree at $x \in \lbrace S_{u_t}, S_{u_t+1},\ldots,S_{\ell_{t+1}-1} \rbrace$. For the linear function, the derivative is a constant of $-1$, and the continuous extension of $f_j(x)$ has derivative of $-1$ at each end of its intervals. Thus, the extension of $\tilde f_t(x)$ is differentiable and the one-sided derivatives are $-1$ at the end of the interval on which it is defined.

To check local Lipschitz continuity of $\tilde f_t(x)$, we note that each component of $\tilde f_t(x)$ is locally Lipschitz continuous in its domain either because it is a linear function or by \cref{result-divergence-component-properties}. Hence, we need to only check that local Lipschitz continuity holds for each $x \in \lbrace S_{u_t}, S_{u_t+1},\ldots,S_{\ell_{t+1}-1} \rbrace$. For $j \in \lbrace u_t,\ldots,\ell_{t+1} - 1\rbrace$, the derivative of $\tilde f_t(x)$ is $-1$ in $(S_j-m_j/32, S_j+m_j/32)$. Hence, the derivative is locally Lipschitz continuous at the stated values of $x$.
\end{proof}

\begin{proposition} 
The function $F: \mathbb{R} \to \mathbb{R}$ as defined in \cref{eqn-divergence-objective-example} is continuous and differentiable on its domain; it is lower bounded; its derivative is locally Lipschitz continuous; $F(S_{\ell_t}) = O_{\ell_t}$, $\forall t \in \mathbb{N}$; $F(S_{u_t}) = O_{u_t}$ $\forall t \in \mathbb{N}$; and $F$'s derivative is not globally Lipschitz continuous.
\end{proposition}
\begin{proof}
The proof is similar to \cref{result-divergence-objective-properties}. Hence, we will only verify that $F$ is lower bounded. By \cref{result-divergence-component-objective-example}, the component $\tilde f_t(X)$ of $F$ for some $t +1 \in \mathbb{N}$ is bounded from below by some $O_j - \frac{3m_{j}}{32}$ for some choice of $j$. So it is enough for us to show, $\lbrace O_j - \frac{3m_j}{32} \rbrace$ is bounded from below. By construction, $\liminf_j O_j = 0$ and $\lim_{j} m_j = 0$. Hence, $\lbrace O_j - \frac{3m_j}{32} \rbrace$ is bounded from below. Thus, $F$ is bounded from below.
\end{proof}

\subsection{Properties of Gradient Descent on the Objective Function}
We now show that when gradient descent is applied to the constructed problem, the objective function's limit supremum is infinite and limit infimum is zero, all while the gradient function remains bounded away from $0$.

\begin{proposition}
Let $\lbrace m_k : k +1 \in \mathbb{N} \rbrace$ be any positive sequence such that $\sum_{k} m_k$ diverges and $m_k \to 0$. Define $F: \mathbb{R} \to \mathbb{R}$ as in \cref{eqn-divergence-objective-example}. Suppose $x_0 = 0$ and let $\lbrace x_k : k \in \mathbb{N} \rbrace$ be generated according to \cref{eqn-recursion-iterate} with $M_k = m_k I$ for all $k+1 \in \mathbb{N}$. Then, $\lbrace M_k \rbrace$ satisfies \cref{property-spd,property-diverge,property-diminishing}. Moreover, (a) $\lim_k x_k = \infty$; (b) $\limsup_k F(x_k) = \infty$; (c) $\liminf_k F(x_k) = 0$; and (d) $\lim_k |\dot F(x_k)| = -1$.
\end{proposition}
\begin{proof}
We first show, $x_k = S_k$ for all $k \in \mathbb{N}$. $0 = x_0 = S_0$. Suppose the claim is true up to $k \in \mathbb{N}$. Then, $\exists t+1 \in \mathbb{N}$ such that $\dot F(x_k) = \dot {\tilde f_t}(x_k)$. Using \cref{result-divergence-component-objective-example} or properties of a linear function, $\dot F(x_k) = \dot F(S_k) = \dot{\tilde f}_t(S_k) = -1$. Therefore,
\begin{equation}
x_{k+1} = x_k - M_k \dot F(x_k) = S_k - m_k \dot{\tilde f}_t(S_k) = S_k + m_k = S_{k+1}.
\end{equation}
Thus, as $k \to \infty$, the iterates diverge and $\dot F(x_k) = -1$ for all $k+1 \in \mathbb{N}$. Now, $F(x_k) = F(S_k) = O_k$ for every $k+1 \in \mathbb{N}$. By properties of $\lbrace O_k \rbrace$, the limit supremum and limit infimum of this sequence is $\infty$ and $0$, respectively. The result follows. 
\end{proof}

We stress that the choice of the limit supremum and limit infimum can be readily modified by choosing a different definition for $\lbrace \ell_t \rbrace$ and $\lbrace u_t \rbrace$. Hence, the limit supremum can be made to be finite and even agree with the limit infimum. Moreover, the limit infimum can be set larger than $0$.

\section{Divergence Regime: Objective Function Diverges, Gradient Function Converges to Zero} \label{subsection-divergence-gradient-zero}
Here, we construct an objective function that is bounded below and has locally Lipschitz continuous gradients. Importantly, when we apply gradient descent with diminishing step sizes to this objective function, the iterates of the procedure diverge, the objective function evaluated at the iterates will diverge, and the gradient function will converge to zero. This objective function will be constructed in a similar fashion to our other divergence regime examples.

\subsection{Construction of the Objective Function}

Let $\lbrace m_{k} : k+1 \rbrace$ be a positive sequence such that $\sum_{k} m_k$ diverges and $m_k \to 0$. Let $S_0 = 0$ and $S_{k+1} = \sum_{j=0}^k m_j$. Then, $\sum_{k} \frac{m_k}{S_{k+1}}$ diverges. Let $K = \min\lbrace k > 0 : S_k \geq 1 \rbrace$ and define
\begin{equation}
T_k = \begin{cases}
S_k & k=0,\ldots,K, \\
T_K + \sum_{j=K}^{k} \frac{m_j}{S_{j+1}} & k > K.
\end{cases}
\end{equation}
Moreover, define
\begin{equation}
d_k = \begin{cases}
1 & k = 0,\ldots,K, \\
\frac{1}{S_{k+1}} & k > K. 
\end{cases}
\end{equation}

Finally, let
\begin{equation} \label{eqn-divergence-gradient-convergence}
F(x) = \begin{cases}
-x & x \leq 0,  \\
f_j(x) & x \in (T_j, T_{j+1}], ~j+1 \in \mathbb{N},
\end{cases}
\end{equation}
where $f_0(x),\ldots,f_{K-1}(x)$ are identical to \cref{eqn-divergence-component-general}; and, letting $x' = x - T_j$,
\begin{small}
\begin{equation} \label{eqn-divergence-gradient-convergence-component}
\begin{aligned}
&f_j(x) \\
&= \begin{cases}
-d_jx' + f_{j-1}(T_j) & x' \in \left[ 0, \frac{(2-d_j) m_j}{16 S_{j+1}} \right) \\
\frac{8S_{j+1}}{m_j}\left( x' - \frac{m_j}{8 S_{j+1}} \right)^2 - \frac{m_j}{S_{j+1}}\left( \frac{-d_j^2 +4d_j}{32} \right) + f_{j-1}(T_j) & x' \in \left[ \frac{(2-d_j) m_j}{16 S_{j+1}}, \frac{3m_j}{16 S_{j+1}} \right) \\
\frac{-5m_j}{16 S_{j+1}} \exp\left( \frac{5/16}{S_{j+1}x'/m_j - 1/2} + 1 \right) + \frac{m_j}{S_{j+1}} \left( \frac{11 + d_j^2 - 4d_j}{32} \right) + f_{j-1}(T_j) & x' \in \left[\frac{3m_j}{16 S_{j+1}}, \frac{m_j}{2S_{j+1}} \right) \\
\frac{m_j}{S_{j+1}}\left(\frac{11 + d_j^2 - 4d_j}{32}  \right) + f_{j-1}(T_j) & x' = \frac{m_j}{2S_{j+1}} \\
\frac{5m_j}{16 S_{j+1}} \exp\left(\frac{-5/16}{S_{j+1}x'/m_j - 1/2} +1 \right) + \frac{m_j}{S_{j+1}}\left( \frac{11 + d_j^2 - 4d_j}{32}  \right) + f_{j-1}(T_j) & x' \in \left( \frac{m_j}{2S_{j+1}}, \frac{13m_j}{16 S_{j+1}} \right) \\
\frac{-8S_{j+1}}{m_j} \left( x' - \frac{7m_j}{8S_{j+1}} \right)^2 + \frac{m_j}{S_{j+1}} \left(\frac{22 + d_j^2 - 4d_j}{32}  \right) + f_{j-1}(T_j) & x' \in \left[ \frac{13 m_j}{16 S_{j+1}}, \frac{(d_{j+1} + 14)m_j}{16 S_{j+1}} \right) \\
- d_{j+1} x' + \frac{m_j}{S_{j+1}} \left(\frac{22 + d_j^2 + d_{j+1}^2 - 4d_j + 28d_{j+1}}{32}  \right) + f_{j-1}(T_j) & x' \in \left[ \frac{(d_{j+1} + 14)m_j}{16 S_{j+1}}, \frac{m_j}{S_{j+1}} \right],
\end{cases}
\end{aligned}
\end{equation}
\end{small}
for $j \geq K$.

\subsection{Properties of the Objective Function}
Here, we verify, \cref{eqn-divergence-gradient-convergence} satisfies \cref{assumption-lower-bound,assumption-gradient}.
We begin by studying the properties of $f_{j}(x)$ for $j \geq K$. Note, we already know the properties of $f_j(x)$ for $j < K$ by \cref{result-divergence-component-properties}.

\begin{proposition} \label{result-divergence-gradient-convergence-component}
Let $j > K$. The continuous extension of $f_j: (T_j, T_{j+1}] \to \mathbb{R}$, \cref{eqn-divergence-gradient-convergence-component}, to $[T_j,T_{j+1}]$ is continuous on its domain; bounded from below by $f_{j-1}(T_j) - m_j/(8S_{j+1})$; differentiable on its domain with $\dot f_j(T_j) = -d_j$ and $\dot f_j(T_{j+1}) = -d_{j+1}$; its derivative is locally Lipschitz continuous; and $f_j(T_{j+1}) \geq f_{j-1}(T_j) + 7m_j/(16 S_{j+1})$.
\end{proposition}
\begin{proof}
The proof of this result is similar to that of \cref{result-divergence-component-properties}. Hence, we only produce the values of $f_j(x)$ and $\dot f_j(x)$ at key points.
\begin{remunerate}
\item At $x = T_j$, $f_j(T_j) = f_{j-1}(T_j)$. $\dot f_j(T_j) = -d_{j}$.
\item At $x = (2-d_j)m_j/(16 S_{j+1}) + T_j$, 
\begin{equation}
f_j(x) = \frac{m_j}{S_{j+1}} \left( \frac{d_j^2 - 2d_j}{16}  \right) + f_{j-1}(T_j),
\end{equation}
and $\dot f_j(x) = -d_j$.
\item At $x = T_j + 3m_j/ (16 S_{j+1})$, 
\begin{equation}
f_j(x) = \frac{m_j}{S_{j+1}} \left( \frac{1 + d_j^2 - 4d_j}{32} \right) + f_{j-1}(T_j),
\end{equation}
and $\dot f_j(x) = 1$.
\item At $x = T_j + m_j/(2 S_{j+1})$,
\begin{equation}
f_j(x) = \frac{m_j}{S_{j+1}} \left( \frac{11 + d_j^2 - 4d_j}{32} \right) + f_{j-1}(T_j),
\end{equation}
and $\dot f_j(x) = 0$.
\begin{equation}
f_j(x) 
\end{equation}
\item At $x = T_j + 13m_j/(16 S_{j+1})$,
\begin{equation}
f_j(x) = \frac{m_j}{S_{j+1}} \left( \frac{21 + d_j^2 - 4d_j}{32}  \right) + f_{j-1}(T_j),
\end{equation}
and $\dot f_j(x) = 1$.
\item At $x = T_j + (d_{j+1} + 14)m_j/(16 S_j)$,
\begin{equation}
f_j(x) = \frac{m_j}{S_{j+1}} \left( \frac{22 + d_j^2 - d_{j+1}^2 - 4d_j}{32} \right) + f_{j-1}(T_j),
\end{equation}
and $\dot f_j(x) = -d_{j+1}$.
\item At $x = T_j + m_j/S_{j+1}$,
\begin{equation}
f_j(x) = \frac{m_j}{S_{j+1}} \left( \frac{22 + d_j^2 + d_{j+1}^2 - 4d_j - 4d_{j+1}}{32}  \right) + f_{j-1}(T_j),
\end{equation}
and $\dot f_j(x) = -d_{j+1}$.
\end{remunerate}
Note, $f_j(T_{j+1}) = f_j(T_j + m_j/S_{j+1}) \geq (22 - 8) m_j / (32 S_{j+1}) + f_{j-1}(T_j)$.
\end{proof}

\begin{proposition}
The function $F: \mathbb{R} \to \mathbb{R}$ as defined in \cref{eqn-divergence-gradient-convergence} is continuous and differentiable on its domain; it is lower bounded; its derivative is locally Lipschitz continuous; $F(T_j) \geq 7T_j/16$ for $j+1 \in \mathbb{N}$; $\dot F(T_j) = -d_{j}$ for all $j+1 \in \mathbb{N}$; and $F$'s derivative is not globally Lipschitz continuous.
\end{proposition}
\begin{proof}
As the proof of this statement is similar to the other constructions, we only verify the values of the objective and the derivative at $\lbrace T_j \rbrace$. For $j=0$, $F(T_0) = 0$. For $j=1,\ldots,K$, $F(T_j) =  F(S_j) = f_{j-1}(S_j) = S_j/2 \geq 7S_j / 16 = 7T_j / 16$ by \cref{result-divergence-component-properties}. For $j > K$, 
$F(T_j) = f_{j-1}(T_j) \geq f_{j-1}(T_{j-1}) + \frac{7 m_j}{16 S_{j+1}} = F(T_{j-1}) + \frac{7m_j}{16 S_{j+1}}$ by \cref{result-divergence-gradient-convergence-component}. By induction, for $j+1 \in \mathbb{N}$, $F(T_j) \geq 7 T_j / 16$. Similarly, either by \cref{result-divergence-component-properties} or \cref{result-divergence-gradient-convergence-component}, $\dot F(T_j) = \dot f_{j}(T_j) = -d_j$.
\end{proof}

\subsection{Properties of Gradient Descent on the Objective}
We now show that when gradient descent is applied to the constructed problem, the objective function diverges, and the gradient function converges to zero.

\begin{proposition}
Let $\lbrace m_k : k +1 \in \mathbb{N} \rbrace$ be any positive sequence such that $\sum_{k} m_k$ diverges and $m_k \to 0$. Define $F: \mathbb{R} \to \mathbb{R}$ as in \cref{eqn-divergence-gradient-convergence}. Suppose $x_0 = 0$ and let $\lbrace x_k : k \in \mathbb{N} \rbrace$ be generated according to \cref{eqn-recursion-iterate} with $M_k = m_k I$ for all $k+1 \in \mathbb{N}$. Then, $\lbrace M_k \rbrace$ satisfies \cref{property-spd,property-diverge,property-diminishing}. Moreover, (a) $\lim_k x_k = \infty$; (b) $\lim_k F(x_k) = \infty$; and (c) $\lim_k |\dot F(x_k)| = 0$.
\end{proposition}
\begin{proof}
We show that $x_k = T_k$ for all $k + 1 \in \mathbb{N}$. For the base case, $x_0 = 0 = T_0$. Suppose $x_k = T_k$ for some $k < K$. Then,
\begin{equation}
x_{k+1} = x_k - M_k \dot F(x_k) = T_k + m_k d_k = T_k + m_k = T_{k+1}.
\end{equation}
This implies that $x_K = T_K$. Now, suppose $x_k = T_k$ for some $k > K$. Then,
\begin{equation}
x_{k+1} = x_k - M_K \dot F(x_k) = T_k + m_k d_k = T_k + \frac{m_k}{S_{k+1}} = T_{k+1}.
\end{equation}
Hence, $F(x_k) = F(T_k) \geq 7T_k / 16$, which diverges to infinity. Moreover, for $k > K$, $|\dot F(x_k)| = |\dot F(T_k)| = d_{k} = 1/S_{k+1}$ which tends to zero.
\end{proof}

\end{document}